%% file: stmca_R2_arxiv.tex
\documentclass[12pt]{article}

\usepackage[T1]{fontenc}
\usepackage[utf8]{inputenc}
\usepackage[nocompress]{cite}
\usepackage[numbers]{natbib}

\usepackage{soul}
\usepackage{cancel}

\usepackage[x11names]{xcolor}
\usepackage{tikz}
\usepackage{pgfplots}

\usepackage{hyperref}
\usepackage{multirow}
\usepackage{hhline}

\usepackage[normalem]{ulem}
\usepackage{amsfonts}
\usepackage{amsmath}
\usepackage{amsthm}
\usepackage{amssymb}
\usepackage{bbm}
\usepackage{mathrsfs}
\usepackage{microtype}

\usepackage{changes}
\usepackage{fullpage}

\usepackage{xcolor}
\usepackage{graphicx}
\usepackage{bm}
\usepackage{enumitem}
\usepackage{verbatim}
\usepackage{algpseudocode}
\usepackage{algorithm}

\usepackage{subcaption}
\usepackage{empheq}
\usepackage{csquotes}
\usepackage{mysymbols2}
\usepackage{siunitx}
\usepackage{authblk}

\usepackage{subfiles}

\newcommand{\stkout}[1]{\ifmmode\text{\sout{\ensuremath{#1}}}\else\sout{#1}\fi}

\theoremstyle{plain}
\newtheorem{theorem}{Theorem}[section]
\newtheorem{corollary}[theorem]{Corollary}
\newtheorem{lemma}[theorem]{Lemma}

\newtheorem{proposition}[theorem]{Proposition}

\theoremstyle{definition}

\newtheorem{example}[theorem]{Example}

\newtheorem{condition}{Condition}

\theoremstyle{remark}
\newtheorem{remark}{Remark}

\DeclareMathOperator{\Bernoulli}{Bernoulli}

\begin{document}

\title{General diffusion processes\\
as the limit of time-space Markov chains}

\author[1]{Alexis Anagnostakis\thanks{alexis.anagnostakis@univ-lorraine.fr}}
\author[1]{Antoine Lejay\thanks{antoine.lejay@univ-lorraine.fr}}
\author[1]{Denis Villemonais\thanks{denis.villemonais@university.edu}}

\affil[1]{Université de Lorraine, CNRS, Inria, IECL, F-54000 Nancy, France}

\renewcommand\Authands{ and }
\date{
}
\maketitle

\subfile{stmca_R2_txt.tex}

	\bibliography{scigenbibfile}
	\bibliographystyle{plainnatnourl}



\end{document}

%% file: stmca_R2_txt.tex
\begin{abstract}
	We prove the convergence of the law of grid-valued random walks, which can be seen as time-space Markov chains,  to the law of a general diffusion process. 
	This includes processes with sticky features, reflecting or absorbing boundaries and skew behavior.  
	We prove that the convergence occurs at any rate strictly inferior to $(1/4) \wedge (1/p)$ in terms of the maximum cell size of the grid, for any $p$-Wasserstein distance. 
	We also show that it is possible to achieve any rate strictly inferior to $(1/2) \wedge (2/p)$ if the grid is adapted to the speed measure of the diffusion, which is optimal for $p\le 4 $. 
	This result allows us to set up asymptotically optimal approximation schemes for general diffusion processes.
	Last, we experiment numerically on diffusions that exhibit various features.
\end{abstract}






\section{Introduction}
In the diffusion process literature, the most well-studied and straightforward
way to approximate diffusion processes is the Euler scheme.  While such
approximations works well for non-degenerate stochastic differential equations,
this is not the case for more general diffusion processes \cite{Hutz2015}. 
The Euler scheme is also not well-defined for processes that
exhibit sticky features, skew behavior \cite{Frik2018} or slowly
reflecting boundaries.

First studied by Feller \cite{Fel52,Fel57}, general diffusion processes
are used to replicate a wide range of phenomena like semi-permeable layers 
\cite{Lej04}, principal agent problem dynamics \cite{Zhu2013,PisWes} and
interest rates behavior close to $0$ \cite{Nie2}. For more theoretical
results and applications of general diffusions, we refer the reader to
\cite{Lej06,EbeZim,Bass,EngPes,HajCagArn}. The resurgence of interest in these
processes and the need to have a universal simulation method for diffusions
have motivated the search for new kinds of approximation processes.

Several works aim at overcoming the shortcomings of the Euler scheme and allow
us to approximate the law of more general diffusion processes.  In \cite{Ami},
the author proposes to approximate the sticky Brownian motion with a simple
random walk that stops for a fixed amount of time when it hits $0$.  In
\cite{Nie2} and \cite{MeoLiGon}, Continuous Time Markov Chains are used to
approximate slowly reflected SDE solutions, where the jumping intensities are
computed using approximated discretizations of the infinitesimal generator of
the diffusion.  Another work where such processes are defined is
\cite{FerLeuJenWeg}, where the authors use a Continuous Time Markov Chain to
identify events in genomics evolution.  In \cite{AnkKruUru,AnkKruUru2}, the
authors develop a numerical scheme to approximate diffusions on natural scale
as long as a mild non-explosion condition is satisfied.  They use symmetric
random walks with fixed-time step whose magnitude depends on the average local
behavior of the target process speed measure.  Choosing the step-size this way
allows the approximation process to replicate non-boundary sticky features.

In this paper, we prove the convergence in law of grid-valued random walks to
any one-dimensional general diffusion process at an asymptotically optimal
rate.  General diffusion processes are regular one-dimensional strong Markov
processes with  continuous trajectories, see for instance \cite[Chapter
7.3]{RevYor} where they are defined as linear continuous Markov processes.
This convergence result allows us to set up approximation schemes that, while
make it straightforward to take into account for sticky points, can also be
applied to any diffusion process that satisfy a mild non-explosion condition.
This includes processes with boundary conditions like absorption, reflection or
slow-reflection as well as the skew diffusions such as the Skew Brownian
motion~\cite{Lej06} and its generalizations.  The values taken by the random
walk correspond to values taken by the target process at random times, allowing
us to classify it as an embeddable scheme along with~\cite{AnkKruUru} and~\cite{EtoLej}.  
We prove that for a grid adapted to the speed measure of the
diffusion process, the laws of the random walks converges at any rate strictly
inferior to  $(1/2) \wedge (2/p)$ in terms of the maximum cell size for all $ p$-Wasserstein distances.  
This convergence rate is
optimal for $p\le 4 $ according the Donsker invariance principle, as this is the rate
simple randoms walk converges to the standard Brownian motion~\cite{Don}. 

Besides the asymptotic optimal convergence rate, the usage of such an
approximation process yield several advantages.  Firstly, the static character
of the grid makes involved quantities good candidates for numerical
approximation (see Sections \ref{sec:setup} and \ref{sec:example}).  Moreover,
this scheme makes it straightforward to take into account potential sticky
points of the diffusion.  Finally, its universality is further validated by the
fact that the Donsker invariance principle and \cite{Don,Ami,EtoLej} are all
special cases of it.

We make heavy use of the speed measure/scale function characterization of diffusions: 
from Proposition 3.14 of \cite[Chapter VII]{RevYor}, the law of a one-dimensional 
diffusion process with state-space~$\II$ an open interval of $\mathbb{R}$ is entirely 
determined by an increasing continuous function $s$ and a positive locally finite measure $m$ defined on $\II$. 
The function $s$ called the \emph{scale function}, is defined as the unique\footnote{up to an affine transformation} 
function such that,
\begin{equation*}
\Prob_x(\tau_b<\tau_a) = \frac{s(x)-s(a)}{s(b)-s(a)},
\end{equation*}
for $a<x<b \in \II$ and with $\tau_a $ being the hitting time of $a $. 
The measure $m$, called the \emph{speed measure}, is the unique positive locally finite measure such that,
\begin{equation*}
\Esp_x(\tau_{ab})=\int_{(a,b)}G_{a,b}(x,y)m(\d y),
\end{equation*}
for $a<x<b \in \II$ and where $G_{a,b}(x,y) $ is the \emph{Green function} defined as
\begin{equation}\label{eq_Green_fun}
G_{a,b}(x,y)=
\begin{cases}
	\dfrac{\big(\scale (x)-\scale (a)\big)\big(\scale (b)-\scale (y)\big)}{\scale (b)-\scale (a)},&\text{for }x\le y,\\
	\dfrac{\big(\scale (y)-\scale (a)\big)\big(\scale (b)-\scale (x)\big)}{\scale (b)-\scale (a)},& \text{for }x> y.
\end{cases}
\end{equation}
We can also express the infinitesimal generator of the diffusion $\Lop$ in terms of $m $ and $s $: 
\begin{align*}
\Lop &= \Dm \Ds,\\
\dom(\Lop) &= \{ f \in \Cob (\II,\mathbb{R}) : \Lop f \in \Cob (\II,\mathbb{R})   \} ,
\end{align*}
where $\Dm$ an $\Ds$ are defined by
\begin{align*}
\Dm f(x)  &= \lim_{h \rightarrow 0;h>0} \frac{f(x+h)-f(x)}{\speed (x,x+h]},\\
\Ds f(x)  &= \lim_{h \rightarrow 0;h>0} \frac{f(x+h)-f(x)}{\scale(x+h)-\scale(x)}.
\end{align*}
These results give us analytic formulations of  quantities of the form 
$v_k(x)=\Esp_x(\tauab^k \indicB{{\tau_b < \tau_a}}) $
for $k\in \mathbb{N}_0 $
and allow us to bound these quantities in terms of the size of the interval $(a,b) $. The latter will be particularly useful for 
proving the convergence results of the algorithm.

\bigskip

\textbf{Outline.} In Section~\ref{sec:scheme}, we present the approximation
scheme along with its properties. In Section~\ref{sec:moment}, we give
analytical characterizations of the quantities that determine the law of the
random walk defined by the algorithm, allowing us to implement it.
Section~\ref{sec:times} is dedicated to proving the convergence of  embedding
times. In Section~\ref{sec:convergence}, we prove the main convergence result
in terms of the maximum cell size of the grid.  The case of solution of an SDE
is studied in Section~\ref{sec:setup}.  Section~\ref{sec:example} is dedicated
to numerical experiments.  In Section~\ref{ssec_lawapprox}, we apply the scheme
for several types of diffusions (reflection, stickiness, skew-behavior,
singularity at 0 or $\infty $) and illustrate its convergence.
In Section
\ref{ssec_localtime}, we exhibit the flexibility of STMCAs on a local time
approximation problem.

\section{The Space-Time Markov Chain Approximation and its properties}\label{sec:scheme}

\subsection{The approximation scheme}\label{ssec:scheme}
In this section we define the approximation process for a one-dimensional diffusion process 
on natural scale\footnote{Which means that $s(x)=x $. } with state space $\II$, an open interval of $\IR$.  
The general case is obtained by a change of scale, as detailed in Section~\ref{ssec:general}.

\begin{figure}[h!]
\centering
\begin{tikzpicture}[scale=0.75]
	
	\draw (-5.5,-2.5)-- (5.5,-2.5); 
	\foreach \x in {-5,-4,-3,-2,-1,0,1,2,3,4,5} {
		\draw (\x,-2.2) -- (\x,-2.8) node[below] {$\x$};
	}
	
	\draw (-5.5,-1.5) -- (-4,-1.5);
	\draw[fill=white] (-4,-1.5) circle (0.09);
	
	\draw (-5,-1.0) -- (-3,-1.0);
	\draw[fill=white] (-3,-1) circle (0.09);
	\draw[fill=white] (-5,-1) circle (0.09);
	
	\draw (-4,-0.5) -- (-2,-0.5);
	\draw[fill=white] (-4,-0.5) circle (0.09);
	\draw[fill=white] (-2,-0.5) circle (0.09);

	\draw (-3,0) -- (-1,0);
	\draw[fill=white] (-3,0) circle (0.09);
	\draw[fill=white] (-1,0) circle (0.09);

	\draw (0,0.5) -- (-2,0.5);
	\draw[fill=white] (-2,0.5) circle (0.09);
	\draw[fill=white] (0,0.5) circle (0.09);
	
	\draw (1,1) -- (-1,1);
	\draw[fill=white] (-1,1) circle (0.09);
	\draw[fill=white] (1,1) circle (0.09);
	
	\draw (0,1.5) -- (2,1.5);
	\draw[fill=white] (0,1.5) circle (0.09);
	\draw[fill=white] (2,1.5) circle (0.09);
	
	\draw (1,2) -- (3,2);
	\draw[fill=white] (1,2) circle (0.09);
	\draw[fill=white] (3,2) circle (0.09);
	
	\draw (2,2.5) -- (4,2.5);
	\draw[fill=white] (2,2.5) circle (0.09);
	\draw[fill=white] (4,2.5) circle (0.09);
	
	\draw (3,3) -- (5,3);
	\draw[fill=white] (3,3) circle (0.09);
	\draw[fill=white] (5,3) circle (0.09);
	
	\draw (4,3.5) -- (5.5,3.5);
	\draw[fill=white] (4,3.5) circle (0.09);

\end{tikzpicture}

\caption{\label{fig:covering}
	The covering grid $\bg = \IZ $ of $(-\infty,\infty) $ along with some of its cells $\bC(\bg) $.}
\end{figure}
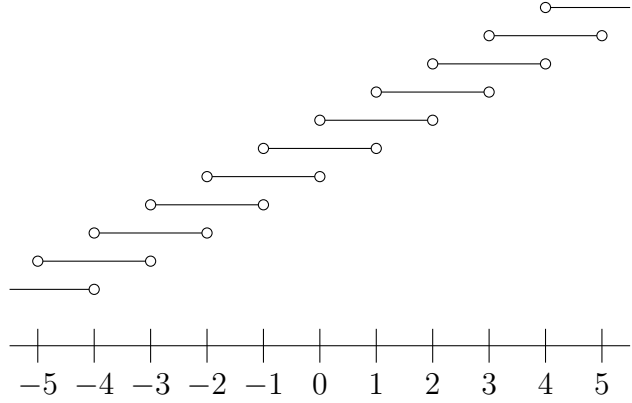

The possible values taken by the approximation process are given as input of the scheme and must form a \emph{covering grid} of $\II $.
We introduce incrementally this notion as follows, which we illustrate by Figure~\ref{fig:covering}.
Let $\II $ be an interval of $\IR $.
\begin{itemize}
\item A \emph{grid} over $\II $ is a countable subset of $\II $ with no accumulation points within $\II $.
\item A \emph{cell} $\bc $ of a grid $\bg $ is an open interval with endpoints in $\bg $ with a single element of $\bg $ in its interior, \textit{i.e.},
\begin{equation*}
	\card\bigbraces{\bg \cap \bc} = 1.
\end{equation*}
We denote with $\bC(\bg) $ the set of all cells of the grid $\bg $.
\item  We call $x\in \II$ the \emph{center} of the cell $\bc \in \bC(\bg) $ iff
\begin{equation*}
	\bc \cap \bg=\{x\}.
\end{equation*}
\item Finally, a \emph{covering grid} of $\II $ is a grid $\bg $ such that 
$\II = \bigcup_{\bc\in \bC(\bg)}\bc $.
\end{itemize}

Examples of covering grids of $(0,\infty) $ and $(-\infty,\infty) $ are $\{1/n; n\in \IN\} $ and $\IZ $ respectively.

For any covering grid $\bg $ of $\II $ and diffusion process $X$ with state-space $\II $, defined through $(s,m)$, let  $|\bg| $, $\cxnorm{X}{\bg} $ be the  grid metrics:
\begin{equation}
\label{eq_grid_metrics_def}
|\bg| = \sup_{\bc \in \bC(\bg)}|\bc| ,\qquad \cxnorm{X}{\bg} = \sup_{\bc \in \bC(\bg)} \{s(\bc)m(\bc)\},
\end{equation}
where $|\bc| = (b-a) $ and $s(\bc) = s(b)-s(a)$ with $a$ and $b$ being the endpoints of $\bc $.
The convergence results will be expressed in the latter metric.

Let $X$ be a diffusion process on natural scale with state-space $\II $ and speed measure $m$.
Let $X$ be defined on the family of filtered probability spaces $(\Omega,\mathcal{F},(\mathcal{F}_t)_{t\ge 0}, (\Prob_x)_{x\in \II})$, where  $\Prob_x$ is the law of $X$ such that, for any $x\in \II $, $\Prob_x \{X_0 = x \}=1$.
For any covering grid $\bg$ of $\II $, we define the process $\widetilde{X}^{\bg} = (\widetilde{X}^{\bg}_{t})_{t\ge 0} $ as the asymmetric random walk 
with:
\begin{itemize}[noitemsep,topsep=-\parskip]
\item  state-space $\bg $,
\item   initial distribution equal to the distribution of $X $ the first time it touches the grid,
\item the same transition probabilities as $X$ over $\bg $,
\item conditional transition times that match the conditional expected transition times of $X $ over $\bg$.
\end{itemize}

Thus, under $\Prob_x $, if $a $ and $b$ are respectively the closest lower and upper elements to $x$ of~$\bg$, 
\begin{equation}\label{eq_approximation_process_1}
\wX^{\bg}_0 = 
\begin{cases}
	a, & \text{with probability } \Prob_x(\tau_b < \tau_a),\\
	b, & \text{with probability } \Prob_x(\tau_a < \tau_b) = 1 - \Prob_x(\tau_b < \tau_a),\\
\end{cases}
\end{equation}
where  $\tau_a := \inf\{t>0: X_t = a\}$.
For the rest of the trajectory, we define $\tauab := \tau_a \wedge \tau_b $ and $\dprocess{\Trt{n}{\bg}}$ as the consecutive jumping times of  $\process{\widetilde{X}^{\bg}_{t} }$.
Then, for all  $k\in \mathbb{N}_0 $ and $a<x<b $ adjacent points of $\bg $,
\begin{equation}
\label{eq_def_probtk}
\Prob\big(\wX^{\bg}_{\Trt{k+1}{\bg}} 
= b \big| \wX^{\bg}_{\Trt{k}{\bg}}  
= x\big) = \Prob_x(\tau_b < \tau_a),
\end{equation}
and 
\begin{equation}
\label{eq_def_esptk}
\Trt{k+1}{\bg} - \Trt{k}{\bg} = 
\begin{cases}
	\Esp_x(\tauab | \tau_b < \tau_a), & \text{on } \{ \wX^{\bg}_{\Trt{k+1}{\bg}}  = b  \}  \cap \{  \wX^{\bg}_{\Trt{k}{\bg}}  = x \},\\
	\Esp_x(\tauab | \tau_a < \tau_b), &\text{on } \{ \wX^{\bg}_{\Trt{k+1}{\bg}}  = a  \}  \cap \{  \wX^{\bg}_{\Trt{k}{\bg}}  = x \}. 
\end{cases}
\end{equation}

As proved in Section~\ref{sec:moment}, the quantities that appear on the right hand side of~\eqref{eq_def_probtk} and~\eqref{eq_def_esptk} are explicit functionals of the speed measure $m $.

Let $\bc_x $ be the cell of the grid $\bg $ containing $x$, \textit{i.e.}, $\bc_x \in \mathcal{C}(\bg) $ and $x = \bc_x \cap \bg $. 
From \eqref{eq_def_probtk}, \eqref{eq_def_esptk} and Bayes' rule, if $\bc_x = (a,b) $, both $\Prob\big(\wX^{\bg}_{\Trt{k+1}{\bg}} = b  \big| \wX^{\bg}_{\Trt{k}{\bg}}  = x\big) $ and $\Trt{k+1}{\bg} - \Trt{k}{\bg}$ only depend  on $x$. 
Thus, if we know the quantities
\begin{equation}\label{algo_quantities}
\begin{split}
	p^{+}[x,(a,b)] =\Prob_x(\tau_b <\tau_a),&\qquad T^{+}[x,(a,b)] = \Esp_x\big[ \tau_{ab} \big| \tau_b< \tau_a \big], \\
	p^{-}[x,(a,b)] =\Prob_x(\tau_a <\tau_b),&\qquad T^{-}[x,(a,b)] = \Esp_x\big[ \tau_{ab} \big| \tau_a< \tau_b \big],\\
\end{split}
\end{equation}
for any adjacent $a<x<b $ in $\bg $,
we can simulate the random walk using Algorithm~\ref{algorithm_3}.
We discuss in Section~\ref{sec:moment} on how to compute the quantities in \eqref{algo_quantities}.
Practical examples are given in Sections~\ref{sec:setup} and~\ref{sec:example}.

This algorithm has been first introduced in \cite{EtoLej} in the situation of SDE solutions with measurable coefficients, where the speed measure of the process satisfies 
\begin{equation}
\label{cond_eto_lej}
c \vd x \le m(\d x) \le C \vd x.
\end{equation}
Our main contribution is that we allow non-elliptic speed measures with vanishing and unbounded density (as in e.g. the Bessel process case, see Section~\ref{ssec_lawapprox}), speed measures with singular part (as in e.g. the sticky Brownian motion setting, see Sections~\ref{ssec_lawapprox} and~\ref{ssec_localtime}), scale functions not in $C^{1} $ (as in e.g. the skew Bessel process, see Section~\ref{ssec_lawapprox}) and non-trivial boundary behaviors (see Section~\ref{ssec:boundary}). 
The probabilistic arguments we use to prove our results allow for greater flexibility, while the proofs of \cite{EtoLej} are based on elliptic PDE theory.
This allows us to handle degenerate diffusions and to perform grid tuning and achieve higher orders of convergence (see Section~\ref{ssec:tuning}).

\begin{remark}
We observe that the process $\widetilde{X}^{\bg} $ is not a one-dimensional Markov chain.
It is though a Markov chain in space and time since 
the joint law of the next position of $\widetilde{X}^{\bg} $ on the grid and the next transition time are both determined by the current position on the grid. 
Hence the terms: space-time Markov chains and Space Time Markov Chain Approximation (STMCA).
\end{remark}

\begin{remark}
In the case of sticky diffusions, where the speed measure $m$ has the form $m(\d x) = m_c(\d x) + \rho \delta_{0}(\d x)  $, the transition probabilities and transition times \eqref{algo_quantities} can be directly inferred from the ones of the diffusion without the sticky term. Indeed, Proposition~\ref{pro_vk_equiv} yields
\begin{equation}\label{eq_cexp_time_sticky_sing}
	\Esp_x\big[\tauab \indicB{\tau_b < \tau_a} \big] = \int_{(a,b)}G_{a,b}(x,\zeta) v_0(\zeta) m_c(\d\zeta)
	+ \rho \,G_{a,b}(x,0) v_0(0),
\end{equation}
where $v_0(x) = \Prob_x (\tau_b < \tau_a) $.
\end{remark}

\begin{algorithm}[h]
\caption{\label{algorithm_3} Space-Time Markov Chain Approximation (STMCA) Algorithm}
\begin{algorithmic}
	\State {\bf Input: }$x$ initial value, $T$ time horizon, $\bg = \{x_j\}_{j\in J}$ grid on $\II$
	
	\State {\bf Output: } $(\widehat{X}[k])_k$ consecutive values taken by the approximation process, $(\hat{t}[k])_k$ consecutive transition times
	\State
	
	\State {\bf Initialization: }
	\State  $\hat{t}[0]=0$, $n=0$
	\State $j = \argmin_{i\in J}\{|x_i-x| \} $ 
	\State $U \sim \Bernoulli\big(p^{+}[x,(x_{j},x_{j+1})]\big) $\;
	\If{$U == 1$}
	\State	$j = j +1 $
	\EndIf
	\State $\widehat{X}[0]= x_j$

	\State
	\State {\bf Main loop: }
	\While{$\hat{t}[n] < T$}
	\State	$U \sim \Bernoulli\big(p^{+}[x_j,(x_{j-1},x_{j+1})]\big) $
	\If {$U==1$}
	\State		$j = j+1$
	\State		$\hat{t}[n+1]=\hat{t}[n] + T^{+}[x_{j},(x_{j-1},x_{j+1})] $\;
	\Else
	\State		$j= j-1$
	\State		$\hat{t}[n+1]=\hat{t}[n] + T^{-}[x_{j},(x_{j-1},x_{j+1})]$
	\EndIf
	\State	$\widehat{X}[n+1]=x_{j}$
	\State	$n=n+1$
	\EndWhile
\end{algorithmic}
\end{algorithm}

\bigskip

\noindent
For the convergence, we make the further assumption that the diffusion process $\process{X_t}$ satisfies the following non-explosion condition: 
there exists a $k_1 > 0 $  such that the speed measure of the diffusion process satisfies 
\begin{equation}
\label{eq_condition1}
m(\d x) \ge k_1 \frac{1}{1+ x^2}\vd x, 
\end{equation}
for all $x\in \II $.
Practically, this means that the process does not move faster than a log-normal process for large values of $X_t$.
We may now express the convergence result in terms of the step-size of the grid in terms of $p$-Wasserstein distances.
In the following result, the $p$-Wasserstein distance $\mathcal{W}_p$ between two laws $\mu$ and $\nu$ of processes with càdlàg\footnote{This stands for \textquote{continue à droite avec une limite à gauche}, that is right-continuous with left-limit.}
paths is defined as
\begin{equation*}
\label{eq_wass_def}
\mathcal{W}_p\big[\mu,\nu \big] = \inf_{(\zeta,\xi) \sim \Gamma(\mu,\nu)} \Big\| \|\zeta - \xi \|_{\infty}
\Big\|_{L^p},
\end{equation*}
where by $\Gamma(\mu,\nu) $ we denote the collection of all measures with marginals $\mu $ and $\nu $.

\begin{theorem}\label{thm_main}
Let $X$ be a diffusion process with state-space $\II $ an interval of $\IR $, on natural scale, whose speed measure satisfies Condition~\eqref{eq_condition1} for some constant $k_1  > 0$. 
Let $\bg$ be a covering grid of $\II $. 
Then, for all $ p\ge 1$, $\delta \in (0,\frac{1}{4} \wedge \frac{1}{p}) $, $T>0 $ and $x\in \II$ there exists a constant~$C>0$ such that
\begin{equation}\label{eq_thm_main}
	\mathcal{W}_p\Big[\Law\big ((\widetilde{X}^{\bg}_{t})_{t\in [0,T]} \big) ,\Law\big ((X_{t})_{t\in [0,T]} \big)\Big]
	\le
	C   \xcnorm{X}{\bg}^{\delta},
\end{equation}
where $\xcnorm{X}{\bg} =  \sup_{\bc \in \bC(\bg)} \cubraces{|\bc| \speed \braces{\bc} }$. 
\end{theorem}

\begin{remark}
In the case where $m(\d x)\ge k_1 \vd x $, the constant\footnote{This results in the bound of Theorem \ref{thm_main} not depending on the starting point of the diffusion.} $C>0$ in Theorem \ref{thm_main} does not depend on the starting point of the diffusion. 
\end{remark}

\begin{remark}
If $X$ is a diffusion process on natural scale such that \eqref{cond_eto_lej} holds,
the bound in~\eqref{eq_thm_main} can be replaced by $C|\bg|^{2 \delta} $.
\end{remark}

The convergence of the Wasserstein distances implies the convergence in law \cite[p.~109]{Villani2009}.

\begin{corollary}
For all $T>0 $, the processes $(\widetilde{X}^{\bg}_t)_{t\in [0,T]} $ converges in law to $(X_t)_{t\in [0,T]} $ in the Skorokhod space $D([0,T],\II) $ as $\xcnorm{X}{\bg}\longrightarrow 0$.
\end{corollary}

In Section~\ref{ssec:markov}, we observe that, by suitable choice of the
probability space, the values taken by the approximation process correspond to
values taken by $X$.  We call this class of approximation schemes embeddable
schemes (other embeddable schemes are the ones developed
in~\cite{AnkKruUru,EtoLej}).  Proving the convergence of an embeddable scheme
usually involves: embedding the approximation process in the trajectory of the
target diffusion process and controlling the embedding times, controlling the
speed at which the process moves.

\subsection{Convergence rate for the general case}\label{ssec:general}

The convergence results  established in the previous section are proven in the
case of a diffusion process on natural scale.  In this section, we show how
more general results can be inferred.  Let $X $ be a diffusion process with
state-space $\II $ an open interval of $\IR $, scale function $s $ and speed
measure $m $.  We assume that
\begin{itemize}
\item $s$ belongs to the Sobolev space $W^{1,1}(\II)$, so from Theorem~8.2 of~\cite{Bre0}, as $s$ is continuous, 
\begin{equation*}
	\scale(x)-\scale(y)= \int_{y}^{x}\scale'(t)\vd t,
\end{equation*} 
for all $y,x$  in $\II $.
\item there exists a $k_1>0$ such that for all $x\in \II $,
\begin{equation}\label{condition8}
	\speed(\d x) \ge k_1 \frac{\scale'(x)}{1+ \bigbraces{\scale(x)}^2}\vd x,
\end{equation}  
\item the inverse of $s$ is $\alpha $-Hölder continuous, \textit{i.e.}, there exists a constant $C>0 $ such that for all $x\ne y \in \II $,
\begin{equation}\label{eq_condition2}
	\frac{|\scale^{-1}(\bar{x}) -\scale^{-1}(\bar{y}) |}{|\bar{x}-\bar{y}|^{\alpha}}\le C.
\end{equation}	
\end{itemize}
Given a grid $\bg $, we consider the random walk $\wX^{\bg}_{t} $ defined by Algorithm~\ref{algorithm_3}, where the transition probabilities and transition times in~\eqref{algo_quantities} can be computed using the formulas derived in Section~\ref{sec:moment}. We obtain the following corollary of Theorem \ref{thm_main}.

\begin{corollary}
\label{cor_main}
Let $X$ be a diffusion process with scale function and speed measure satisfying the above conditions. Let also $\bg$ be a covering grid over the state-space $\II$ of $X$. Then, for all $ p\ge 1$, $\delta \in (0,\frac{1}{4}\wedge \frac{1}{p}) $, $T>0 $ and $x\in \II$ there exists positive constants $C_1 $, and $C_2 $ such that
\begin{equation*}
	\label{eq_corr_non_natural_thm_main}
	\mathcal{W}_p\Big[\Law\big ((\widetilde{X}^{\bg}_{t})_{t\in [0,T)} \big) ,\Law\big ((X_{t})_{t\in [0,T)} \big)\Big]
	\le
	C_1 \xcnorm{X}{\bg}^{\delta},  
\end{equation*}
where $\cxnorm{X}{\bg}$ is defined in  \eqref{eq_grid_metrics_def}.
\end{corollary}

\begin{proof}
We define the proxy process $Y = \process{s(X_t)}$  which has state-space
$s(\II) $, scale function $s_Y(x)=x $ and speed measure $m_Y(\d x) = m
\circ s^{-1}(\d x) $.  From condition~\eqref{condition8} and a change of
variables, we get that $\process{Y_t}$ satisfies
condition~\eqref{eq_condition1} for the same constant $k_1$.  We also
define $\widetilde{Y}^{s(\bg)}_{t}$ evolving according to Algorithm~1, with
covering grid $s(\bg) = \{s(x);x\in \bg \} $. It can be defined on the
canonical space of $X$ so that $s(\wX^{\bg}_{t}) =
\widetilde{Y}^{s(\bg)}_{t} $ almost surely. Thus,
Condition~\eqref{eq_condition2} implies that
\begin{equation*}
	\label{eq_Holder_continuity_impl2}
	|X_t - \widetilde{X}_{t}^{\bg}| \le C|Y_t - \widetilde{Y}_{t}^{\scale(\bg)} |^{\alpha }.
\end{equation*}
Along with the fact that
\begin{equation*}
	\label{X_grid_norm}
	\xcnorm{Y}{\bg} = \sup_{\bc \in \bC(\bg)} \{ |s(\bc)| m_Y(s(\bc))  \},
\end{equation*}
Theorem~\ref{thm_main} implies Corollary~\ref{cor_main}.
\end{proof}

\subsection{Grid tuning}\label{ssec:tuning}

We observe that for all $\epsilon>0 $, in the case of a Brownian motion,
Theorem~\ref{thm_main} yields a convergence rate of
$\Ord(|\bg|^{\frac{1}{2}-\epsilon})$ as $|\bg|\lra 0 $, which is optimal from
Donsker's invariance principle \cite{Don}.  We would like to have this
result for all diffusion processes, but the following example illustrates that
this is not the case.  We then show how we can remediate to this by using a
custom grid and extrapolate this method to the general case via
Corollary~\ref{cor_grid_tuning}.

\begin{example}
\label{ex:sticky}
Let $X $ be the diffusion process with state-space $\IR $, defined through $s$ and $m$ with
\begin{equation*}
	s(x) = x, \qquad m(\d x) = 2\vd x + \rho \delta_0(\d x).
\end{equation*}
This process is called the \emph{sticky Brownian motion} and is the \textquote{most elementary} sticky diffusion process. As such, it spends a positive amount of time at $0$ and the Euler scheme is known to not be well defined for these processes.
We observe that, for any covering grid $\bg $ of $\IR $ and $\bc \in \bC(\bg) $ with $\bc \ne \bc_0 $ and $\bc_0\in \bC(\bg) $ being the cell containing $0$,
\begin{align}
	\notag
	2|\bc|^2 &= m(\bc)|\bc| \le \sup_{\bc\in C(\bg)} {m(\bc)|\bc|},\\
	\label{eq_GT_grid_at_sticky_point}
	\rho|\bc_0| + 2|\bc_0|^2 &= m(\bc_0)|\bc_0| \leq \sup_{\bc\in C(\bg)} {m(\bc)|\bc|}.
\end{align}
From \eqref{eq_GT_grid_at_sticky_point}, for a uniform grid of step size h,
\begin{equation*}
	2 h \rho  < \xcnorm{X}{\bg}.
\end{equation*}
Thus, for any $\epsilon >0 $, the convergence rate given by Theorem \ref{thm_main} is $\Ord(h^{\frac{1}{4}-\epsilon}) $ as $|\bg|\lra 0 $.
This means that there are functionals of the trajectory for which the convergence rate is much slower for this process in comparison with a standard Brownian motion. In order to remediate to this, we propose a preliminary step to the approximation scheme that involves finding a grid that is \textquote{adapted}  to the speed measure of the process. In the case of the Brownian motion with a sticky point at $0$, such a grid can be defined as  one that has uniform non-adjacent cells to $0$ of size $h$ and with the cell of center $0$ being $(-h^2/2\rho, h^2/2\rho) $, i.e.

\begin{equation}\label{eq_tuned_grid}
	\bg =\Big\{ \bigcup_{k\in\IZ_{+}} \big\{-\frac{h^2}{2\rho}-k\frac{h}{2}\big\} \Big\} 
	\cup  \big\{0\big\} 
	\cup \Big\{ \bigcup_{k\in\IZ_{+}} \big\{\frac{h^2}{2\rho}+k\frac{h}{2}\big\} \Big\}.
\end{equation}
\noindent
As the approximation process is a random walk, for every $k$ steps it makes, it spends $\mathcal{O}(\sqrt{k})$ steps in the cell containing $0$ (see \cite{Chung1949}). Thus, running the algorithm on either grid yields the same algorithmic complexity, whilte the convergence rate is improved to $\Ord(h^{\frac{1}{2}-\epsilon}) $ for the adapted grid. Numerical examples are given in Section~\ref{example_stqBM}.
\end{example}

The general case is covered by the following Corollary:

\begin{corollary}\label{cor_grid_tuning}
Let $X$ be a diffusion process and $\widetilde{X}^{\bg} $ the approximation process defined by Algorithm \ref{algorithm_3}.
Then, if $\bg $ is a grid such that 
\begin{equation}\label{eq_grid_tuning_condition}
	\xcnorm{X}{\bg} \le C |\bg|^2,
\end{equation}
we can bound the $p$-Wasserstein distance between the laws of $(\widetilde{X}^{\bg}_{t})_{t\in [0,T)}$ and $(X_{t})_{t\in [0,T)}$ in Theorem \ref{thm_main} by $|\bg|^{2\delta} $ instead of $\xcnorm{X}{\bg}^{\delta}$. Thus, for all $\epsilon > 0 $, the law of the random walk converges in any $ p$-Wasserstein distance at the rate $\Ord(|\bg|^{(\frac{1}{2} \wedge \frac{2}{p}) -\epsilon})$ instead of $\Ord(|\bg|^{(\frac{1}{4}\wedge \frac{1}{p} ) - \epsilon})$ as $|\bg|\lra 0$.
\end{corollary}

We now show how one can create grids such that \eqref{eq_grid_tuning_condition} holds in the case of homogeneous SDEs. 
Considering sticky and/or skew points is straightforward.

\begin{itemize}
\item For a process whose speed measure satisfies \eqref{cond_eto_lej} and \eqref{condition8} with one or more skew points, no grid modification is required, 
\item For a process whose speed measure satisfies \eqref{cond_eto_lej} and \eqref{condition8} and has a sticky point at $0 $ of stickiness $ \rho>0$, one needs to consider the points $ \{- h^2/\rho,0,  h^2/\rho\} $ to have a tuned grid,
\item The case of a reflection at a boundary is treated in Section~\ref{ssec:boundary}.
\end{itemize}

Let $(\mu,\sigma) $ be a pair of functions satisfying the following condition:
\begin{condition}\label{cond_sde_condition}
The functions are measurable $\IR \mapsto \IR $ mappings and
the SDE
\begin{equation}\label{eq_sde_sol}
	\d X_t = \mu (X_t)\vd t + \sigma(X_t)\vd B_t,
\end{equation}
has a unique weak solution, where $B$ is a standard Brownian motion.
\end{condition}

Let $X$ be the diffusion that solves \eqref{eq_sde_sol}, $\II $ its state-space, $s$ its scale function and $m$ its speed measure given by \cite[p. 17]{BorSal}
\begin{equation}\label{eq_dm_recap}
s'(x)=e^{-\int_{a}^{y}\frac{2\mu (\zeta)}{\sigma^2 (\zeta)}\vd \zeta},
\quad\text{and}\quad
m(\d x) = \frac{1}{s'(x)}\frac{2}{\sigma^2(x)}\vd x,
\end{equation}
with $s'$ being the right-derivative of $s$.
Then, if $\bg =\{x_j\}_{j\in J} $ is a covering grid such that for a constant $C>0 $ and every $j\in J $
\begin{equation}
\Bigabsbraces{ \bigbraces{s(x_{j+1})-s(x_{j-1})}\int_{x_{j-1}}^{x_{j+1}} \frac{2}{s'(\zeta )\sigma^2 (\zeta)} \vd \zeta} \wedge (x_{j+1}-x_{j-1})^2 \le C h^2 ,
\end{equation}
it satisfies \eqref{eq_grid_tuning_condition} and $|\bg|=C h $.
Thus, from Corollary \ref{cor_grid_tuning}, for any $\epsilon>0 $  using such grids give us a convergence rate of $\Ord(|\bg|^{\frac{1}{2} - \epsilon})$ in Theorem \ref{thm_main}.

Generating such grids numerically can be done choosing a starting point $x_0 $ and adding points $x_j $ to the grid  iteratively as follows: 
given $x_{j-1} $, let $x_j $~be chosen such that:
\begin{equation}\label{eq_grid_tuning_criteria_SDE}
\Bigabsbraces{ \bigbraces{s(y)-s(x_{j-1})}\int_{x_{j-1}}^{y} \frac{2}{s'(\zeta)\sigma^2 (\zeta)} \vd \zeta}    \le h^2/2.
\end{equation}
Then the next element of the grid is defined as
\begin{equation}\label{eq_grid_tuning_reiteration}
x_{j} =  
\begin{cases}
	x^{(0)}_{j} &\quad\text{if } x^{(0)}_j - x_{j-1} \le h, \\
	x_{j-1}+h   &\quad\text{if } x^{(0)}_j - x_{j-1} > h. \\
\end{cases}
\end{equation}
The problem \eqref{eq_grid_tuning_criteria_SDE} can be solved numerically using a fixed point algorithm.
Examples of STMCA simulations using tuned grids computed solving \eqref{eq_grid_tuning_criteria_SDE}-\eqref{eq_grid_tuning_reiteration} are given in Figures \ref{fig:grid_tuning_cir3} and \ref{fig:grid_tuning}.

\subsection{The case of diffusions with boundary conditions}\label{ssec:boundary}

When presenting the results and the structure of the scheme, we considered only processes  where $\II$ is an open set, thus considering diffusion with unreachable boundaries. Our results also adapt to the situation where either $\ell$ and/or $r$ are reachable, and in this case some adjustments are needed, depending on the nature of finite boundaries and on the condition at regular boundaries. In order to keep the presentation simple, we assume that the process is on natural scale and that  $\II=[0,+\infty)$ (the adaptation to $\II=(\ell,r]$ or $\II=[\ell,r]$ or $\II=[\ell,r)$ with $\ell\in\mathbb{R}\cup\{-\infty\}$ and $r\in\mathbb{R}\cup\{+\infty\}$  is straightforward).

It is well known (see for
instance Section~5.11 of It\^o's book~\cite{Ito2006})
that the finite boundary~$0$ can be of four types. Setting, for some
fixed $c>0$,
\begin{align*}
\mathcal{I}=\int\int_{0<y<x<c} m(\d x)\vd y,\quad \mathcal{II}=\int\int_{0<y<x<c} m(\d y)\vd x,
\end{align*}
then
\begin{itemize}[noitemsep,partopsep=-\parskip,topsep=0pt]
\item $0$ is an exit boundary if $\mathcal{I} <\infty$ and
$\mathcal{II}=\infty$,
\item $0$ is a regular boundary if $\mathcal{I}<\infty$ and
$\mathcal{II}<\infty$,
\item $0$ is a natural boundary if $\mathcal{I}=\infty$ and
$\mathcal{II}=\infty$,
\item $0$ is an entrance boundary if $\mathcal{I}=\infty$ and
$\mathcal{II}<\infty$.
\end{itemize}
The entrance type been excluded for a finite boundary of a diffusion process on natural scale, and the natural type been considered in the settings of Theorem~\ref{thm_main},  this leaves us with two possible types for the boundary $0$: exit or regular. If $0$ is an exit boundary, then the diffusion process $X$ is absorbed at the boundary $0$. 
If $0$ is a regular boundary, then the diffusion process can either be absorbed or reflected at $0$.
In these cases, the convergence result of Theorem~\ref{thm_main} can be extended by considering a grid $\bg$ on $\II$ containing $0$ and by adapting the dynamics of $\widetilde{X}^{\bg}$ as follows.
The dynamic of $\widetilde{X}^{\bg}$ is the same as in Algorithm~\ref{algorithm_3}, up to the time when it reaches $0$, then: 
\begin{itemize}[leftmargin=1em]
\item  if $0$ is an absorbing boundary (exit or regular), then the result can be immediately extended by stopping ${\widetilde X}^{\bg}$ when it reaches~$0$;
\item if $0$ is a reflecting regular boundary, then the process  $\widetilde{X}^{\bg}$ jumps from $0$ to  $b:=\min\bg_{\,\setminus\{0\}}$  with probability $1$ and after a time $\int_{[0,b)} (b-\zeta)\,m(\d\zeta)$.
We emphasize that in this configuration, $0$ may be a sticky boundary (\textit{i.e.}, with $m(0)\in(0,+\infty)$).
\end{itemize} 

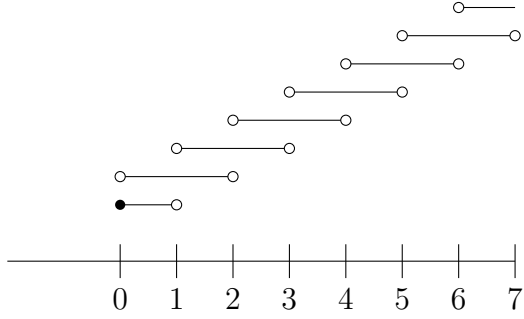
\begin{figure}[h!]
\centering
\begin{tikzpicture}[scale=0.75]
	
	\draw (-2,0)-- (7,0); 
	\foreach \x in {0,1,2,3,4,5,6,7} {
		\draw (\x,0.3) -- (\x,-0.3) node[below] {$\x$};
	}
	\draw (0,1) -- (1,1);
	\fill (0,1) circle (0.09);
	\draw[fill=white] (1,1) circle (0.09);
	
	\draw (0,1.5) -- (2,1.5);
	\draw[fill=white] (0,1.5) circle (0.09);
	\draw[fill=white] (2,1.5) circle (0.09);
	
	\draw (1,2) -- (3,2);
	\draw[fill=white] (1,2) circle (0.09);
	\draw[fill=white] (3,2) circle (0.09);
	
	\draw (2,2.5) -- (4,2.5);
	\draw[fill=white] (2,2.5) circle (0.09);
	\draw[fill=white] (4,2.5) circle (0.09);
	
	\draw (3,3) -- (5,3);
	\draw[fill=white] (3,3) circle (0.09);
	\draw[fill=white] (5,3) circle (0.09);
	
	\draw (4,3.5) -- (6,3.5);
	\draw[fill=white] (4,3.5) circle (0.09);
	\draw[fill=white] (6,3.5) circle (0.09);
	
	\draw (5,4) -- (7,4);
	\draw[fill=white] (5,4) circle (0.09);
	\draw[fill=white] (7,4) circle (0.09);
	
	\draw (6,4.5) -- (7,4.5);
	\draw[fill=white] (6,4.5) circle (0.09);
	
\end{tikzpicture}

\caption{The covering grid $\bg = \IZ^{+} $ of $[0,\infty) $ along with its first cells  $\bC(\bg) $.}
\end{figure}

In both cases of reflection and absorption at $0$, the boundary is attainable. 
Thus, $0$ must be a point of any covering grid of $\II = [0,\infty) $.
This leads to the following adaptation of the notion of grid cells.
The cells $\bC(\bg) $ of  such grid  $\bg $ are the open intervals  for the induced topology of $\IR $ on $\II $ with endpoint in $\bg $.
For example $\bg_0 = \mathbb{Z}_{+} $ is a covering grid of $[0,\infty) $ and $\bC(\bg_0) = \{[0,1), (n-1,n+1)_{ n\in \IN } \} $.

The proof of the convergence in these situations is omitted here, 
since it is a straightforward adaptation of the proof of Theorem~\ref{thm_main}, 
using in particular the fact that, in the case of a reflecting boundary,
\begin{equation*}
\Esp_0\big[\tau_{b}\big] = \int_{[0,b)} (b-\zeta)\,m(\d\zeta).
\end{equation*}
The case of killing boundaries, and in general of a process with non-zero killing measure, 
leads to additional non-trivial difficulties. Devising an algorithm and a similar result 
as Theorem~\ref{thm_main} for such processes remains an active area of research.

\subsection{Markovian embedding}\label{ssec:markov}

The consecutive values of the process $\process{\wX_t} $ defined 
through \eqref{eq_approximation_process_1}-\eqref{eq_def_esptk}
form a Markov chain with, by construction, the same transition probabilities as $ \process{X_t}$ on $\bg $. We define the embedding times of $\process{X_t} $ in $\bg $ as,
\begin{equation}
\label{eq_tauk_def}
\begin{cases}
	\tau^{\bg}_{0}&=0 ,\\
	\tau^{\bg}_{k}&=\inf\big\{t> \tau^{\bg}_{k-1}: X_t \in \bg \setminus \{ X_{\tau^{\bg}_{k-1}}\} \big\}, \qquad  \forall k\ge 1.
\end{cases}
\end{equation}
As both $\wX_{\Trt{n}{\bg}}$ and $X_{\tau^{\bg}_{n}}  $ are both Markov chains with the same transition probabilities with~$\wX_0 $ forced to be equal in law to~$X_{\tau^{\bg}_{1}} $ (see Section~\ref{ssec:scheme}), the following equality in law holds,
\begin{equation}
\Law( \wX_{\Trt{n}{\bg} } ; n \ge 0 )=\Law(X_{\tau^{\bg}_{n}} ; n \ge 1).
\end{equation}
We define $\Kri{t}{\bg} $ as the inverse of $\Trt{n}{\bg}$, \textit{i.e.},
\begin{equation}
\label{eq_random_index_definition}
\Kri{t}{\bg}=\inf \Big\{n\in\mathbb{N}: \sum_{k=1}^{n}\Esp \big[ \tau_{k}^{\bg} - \tau_{k-1}^{\bg} \big| X_{\tau_{k-1}^{\bg}} , X_{\tau_{k}^{\bg}} \big] > t  \Big\}.
\end{equation}	
Thus, we get the following Proposition.

\begin{proposition}
\label{prop_embedding_2}
Let $(X_t)_{t\geq 0}$ be a diffusion process,~$\bg $ a grid defined over its state space~$\II$ and~$\wX^{\bg}_t $
be the approximation process defined in \eqref{eq_def_probtk} and \eqref{eq_def_esptk}. Then, if $(\tau^{\bg}_n)_{n\geq 0}$ are the embedding times of $(X_t)_{t\geq 0}$ in $\bg $, 
the following equality in law holds,
\begin{equation*}
	\Law\big(\wX_t ;t\ge 0 \big) = \Law\big( X_{\tau^{\bg}_{ \Kri{t}{\bg} }} ; t\ge 0 \big),
\end{equation*}
where $\Kri{t}{\bg}$ is the random index defined  in \eqref{eq_random_index_definition}.
\end{proposition}

\section{Moment characterization of conditional exit times}\label{sec:moment}
The law of  the approximation process defined in the previous section was shown to be determined by the transition probabilities $\Prob_x(\tau_b<\tau_a) $ and conditional transition times $\Esp_x(\tau_{ab}| \tau_b < \tau_a) $.  
In this section we show that quantities of the form
$v_k(x)=\Esp_x(\tau_{ab}^k\indicB{ \tau_b < \tau_a})$ yield an integral formulation
with respect to the speed measure of the diffusion and involving the scale function (we do not assume that the diffusion is on natural scale in the present section). 
We also show that this results in them being solutions to  Dirichlet problems
where the differential operator is the infinitesimal generator $\Lop $ of the
diffusion. 
This allows us to simulate such processes via Algorithm~\ref{algorithm_3} and thus to approximate the law of the target diffusion process~$\process{X_t}$.

\noindent
In terms of Algorithm \ref{algorithm_3}, we need to compute for three adjacent points $a,x,b$ 
of the grid the quantities
\begin{equation*}
v_0(x)=\Prob_x(\tau_b<\tau_a),\ v_1(x)=\Esp_x(\tau_{ab}\indicB{\tau_b<\tau_a})
\text{ and }
\overline{v}_1(x)=\Esp_x(\tau_{ab}\indicB{\tau_a<\tau_b}).
\end{equation*}
The quantities of \eqref{algo_quantities} are then 
\begin{equation}
\label{algo_quantities_2}
\begin{aligned}
	p^{+}[x,(a,b)]&=v_0(x), &T^{+}[x,(a,b)] &= \frac{v_1(x)}{v_0(x)},
	\\
	p^{-}[x,(a,b)]&=1-v_0(x),& T^{-}[x,(a,b)] &= \frac{\overline{v}_1(x)}{1-v_0(x)}.
\end{aligned}
\end{equation}

\begin{proposition}\label{prop_v0_DP}
The function $v_0(x)=\Prob_x(\tau_b< \tau_a) $ is solution to the problem with Dirichlet boundary conditions
\begin{equation}\label{eq_prob_P0}
	\begin{cases}
		\Lop u = 0, & x\in (a,b), \\
		u(a) = 0, & \\
		u(b) = 1, & 
	\end{cases}
\end{equation}
where $\Lop$ is the infinitesimal generator of $(X_t)_{t\ge 0}$, which also implies that $v_0 \in \dom(\Lop) $.
\end{proposition}

\begin{proof}
Let $x\in (a,b) $, from the definition of the scale function and the factorization of the infinitesimal generator $\Lop = \Dm \Ds $
\begin{equation*}
	\begin{split}
		\Lop v_0 =  \Dm \Ds \frac{s(\cdot)-s(a)}{s(b)-s(a)}= \Dm \frac{1}{s(b)-s(a)}.
	\end{split}
\end{equation*}	
which equals $0$ as  $m(\d x) $  is a positive measure. As $v_0 $ and $\Lop v_0 =  0 $ are both functions in $C^{0}_{b} $, we deduce that $v_0 \in \dom(\Lop)$ and $\Lop v_0 =  0 $. 
Under $\Prob_b$, the stopping time $\tau_b$ equal $0$ a.s. and the process has a.s. continuous trajectories, hence $\tau_a > 0 $ a.s., \textit{i.e.},
\begin{equation*}
	v_0(b)=\Prob_b(\tau_b<\tau_a)=\Prob_b(0< \tau_a)=1.
\end{equation*}
This, along with the symmetrical argument, allow us to retrieve the boundary conditions of~\eqref{eq_prob_P0}.
\end{proof}

\begin{proposition}
\label{pro_vk_equiv}
For every $k\in \IN $, let $v_k $ be the function defined for every $x\in (a,b) $ by $v_{k}(x)= \Esp_x(\tau_{ab}^k\indicB{ \tau_b < \tau_a})$. Then,
\begin{equation}\label{meq_pro_vk_equiv}
	v_k(x)= k \int_{(a,b)} G_{a,b}(x,\zeta)v_{k-1}(\zeta)m(\d\zeta).
\end{equation}
\end{proposition} 

\begin{proof}
Since $\int_0^{\tau_{ab}}(\tau_{ab}-t)^{k-1}\vd t=\tau_{ab}^k/k$,
\begin{equation*}
	v_k(x) = k \Esp_x \Big[ \indicB{\tau_b< \tau_a} \int_{0}^{\tau_{ab}} (\tau_{ab}-t)^{k-1}  \vd t \Big]= k \Esp_x \Big[ \indicB{\tau_b< \tau_a} \int_{0}^{\infty} \indicB{t \le \tau_{ab}}(\tau_{ab}-t)^{k-1}  \vd t \Big]. 
\end{equation*}
From the Markov property, by conditioning on $\mathcal{F}_t $ and as $\indicB{t \le \tau_{ab}} $ is $\mathcal{F}_t $-measurable, 
\begin{equation*}
	v_k(x)= k \Esp_x \left[  \int_{0}^{\infty} \indicB{t \le \tau_{ab}}  \Esp \big[\indicB{\tau_b< \tau_a} (\tau_{ab}-t)^{k-1}\big| \mathcal{F}_t  \big]  \vd t \right]=k \Esp_x \left[  \int_{0}^{\tau_{ab}}  \Esp_{X_t} \big[\indicB{\tau_b< \tau_a} \tau_{ab}^{k-1}  \big]  \vd t \right].
\end{equation*}
The equality \eqref{meq_pro_vk_equiv} results by applying directly Green's formula.
\end{proof}

\begin{proposition}\label{prop_vk_DP}
The function $v_k(x)=\Esp_x(\tau_{ab}^{k}\indicB{\tau_b<\tau_a}) $ is solution to the problem with Dirichlet boundary conditions
\begin{equation*}\label{eq_prob_Pk}
	\begin{cases}
		\Lop u = -k v_{k-1}, & x\in (a,b), \\
		u(a) = 0, & \\
		u(b) = 0. & 
	\end{cases}
\end{equation*}	
\end{proposition}

\begin{lemma}\label{lem_g_DP}
Let $g(x)= \int_{(a,b)}G_{a,b}(x,y)f(y)m(\d y) $, where $f\in C^{0}_{b}(a,b) $ and $G_{a,b}(x,y)$ is the Green function defined in \eqref{eq_Green_fun}. Then $g\in \dom(\Lop) $ and 
\begin{equation*}
	\Lop g(x) = -f(x), \qquad \forall x\in (a,b).
\end{equation*}
\end{lemma}

\begin{proof}
Let $x\in(a,b) $. Using the $\Dm \Ds $ factorization of $\Lop$ and the dominated convergence theorem we get
\begin{multline*}
	\Lop g(x) 
	\\
	= \Dm \Ds \int_{(a,b)}
	\Big[1 _{y< x}\frac{(s(y)-s(a))(s(b)-s(x))}{s(b)-s(a)} + \indicB{y\ge x}\frac{(s(x)-s(a))(s(b)-s(y))}{s(b)-s(a)} \Big]f(y)m(\d y) \\
	=  -\Dm \int_{(a,x)}v_0(y)f(y)m(\d y) + \Dm \int_{[x,b)}(1-v_0(y)) f(y)m(\d y) \\
	= -v_0(x)f(x)- (1-v_0(x))f(x) = -f(x).
\end{multline*}
The continuity of $g $ is a consequence of Lebesgue's convergence theorem for integrals. Moreover as $G_{a,b}(x,y) $ is bounded by $s(b)-s(a) $, $m(\d x) $ is locally finite and $f $ is bounded, $g$ is also bounded. So, we deduce that $f\in \dom(\Lop) $ and that on $(a,b)$ we have $\Lop g = -f $.
\end{proof}

\begin{proof}[Proof of Proposition \ref{prop_vk_DP}]
From Proposition  \ref{pro_vk_equiv},
\begin{equation*}
	v_k(x)= \int_{(a,b)}G_{a,b}(x,y)kv_{k-1}(y) m(\d y).
\end{equation*}
As $v_0 \in C^{0}_{b}(a,b) $, from Lemma \ref{lem_g_DP} we deduce iteratively that $v_k \in C^{0}_{b}(a,b) $ for all $k\in \mathbb{N} $ and that $\Lop v_k = -k v_{k-1} $ on $(a,b) $. For the boundary conditions, we observe that $\tau_{ab}=\tau_a \wedge \tau_b =0 $ a.s. under $\Prob_a$, so for $k\ge 1 $
\begin{equation*}
	v_k(a)=\Esp_a(\tau_{ab}^{k}\indicB{\tau_b<\tau_a}) \le \Esp_a(\tau_{ab}^{k})=0.
\end{equation*}
With the same argument we show that $v_k(b)\le 0 $ and as they are obviously positive quantities $v_k(a)=v_k(b)=0 $.
\end{proof}

\section{Convergence of the embedding times}\label{sec:times}

In order to prove the convergence of the process $(\widetilde{X}^{\bg}_t)_{t\ge 0} $ we  need to control quantities of the form $\Esp_x \big[ \sup_{t\le T}|t-\tau_{\Kri{T}{\bg}}^{\bg}|^p \big] $, where $\tau_{k}^{\bg} $ are the embedding times of the process $X $ in the grid $\bg $. In this section, we show the existence of such bounds in terms of the metric 
\begin{equation}
\xcnorm{X}{\bg} = \sup_{\bc \in \bC(\bg)} \{s(\bc)m(\bc)\}
\end{equation}
defined in \eqref{eq_grid_metrics_def}. 
If there exists a constant $K> 0 $ such that $m(\bc)\le K $ for all $\bc \in \bC(\bg) $, then, $|\bc|m(\bc) \le K |\bc|  $ and $ \xcnorm{X}{\bg} \le K |\bg|$.
Thus, we can bound the quantities of interest in terms of $|\bg|$ instead of $ \xcnorm{X}{\bg} $.
But, in doing so, we do not track correctly the convergence rates of the approximation process. 
For example, for the standard Brownian motion,  $\xcnorm{B}{\bg} =|\bg|^2 $. Moreover, as shown in Section~\ref{ssec:tuning}, such bounds give us a direct way to adapt the grid to the speed measure in order to accelerate the convergence of the scheme.

\subsection {Bounds on the conditional moments of the exit times}

\begin{lemma}\label{lem_vk_ratio}
Let $v_k(x)=\Esp_x(\tau_{ab}^k \indicB{\tau_b< \tau_a}) $, then for all $k\in \mathbb{N}  $
\begin{equation}\label{eq_vk_ratio_ninf_bound}
	\Big\|\frac{v_k}{v_{k-1}}\Big\|_{\infty} \le k \xcnorm{X}{\bg}.
\end{equation}
\end{lemma}

\begin{proof}
We first observe that
\begin{equation*}
	G_{a,b}(x,y)  \frac{v_{0}(y)}{v_{0}(x)} = 
	\begin{cases}
		\frac{(y-a)(b-x)}{(b-a)}\frac{y-a}{x-a},& x > y,  \\
		\frac{(x-a)(b-y)}{(b-a)}\frac{y-a}{x-a},& x \le y.\\
	\end{cases}
\end{equation*}	
As for $x>y$ the ratio $\frac{y-a}{x-a}<1 $,
\begin{equation}\label{eq_proof_v0_Green_mult_ratios}
	\begin{split}
		G_{a,b}(x,y) \frac{v_{0}(y)}{v_{0}(x)}&\le 
		\begin{cases}
			\frac{(y-a)(b-x)}{(b-a)},& x > y,  \\
			\frac{(y-a)(b-y)}{(b-a)},& x \le y,\\
		\end{cases}
	\end{split}
\end{equation}
which is bounded by $ (b-a) $ in both cases as $x,y\in (a,b) $.
Thus, 
\begin{equation}\label{eq_proof_v0_Green_mult_ratios_integral}
	\int_{a}^{b} G_{a,b}(x,y)v_{0}(y)m(\d y) = 
	v_0 (x)\int_{a}^{b} G_{a,b}(x,y)\frac{v_0 (y)}{v_0 (x)}m(\d y) \le v_0 (x) (b-a) m((a,b)).
\end{equation}
From Proposition~\ref{pro_vk_equiv}, $v_k(x)=k \int_{(a,b)}G_{a,b}(x,y)v_{k-1}(y)m(\d y) $ and
\begin{equation}\label{eq_proof_v0_Green_mult_ratios_integral_B}
	v_k(x)=k! \int_{(a,b)}G_{a,b}(x,x_k) \int_{(a,b)}G_{a,b}(x_k,x_{k-1}) 
	\dots \int_{(a,b)}G_{a,b}(x_2,x_{1})  v_{0}(x_1)m(\d x_1)\dots m(\d x_{k}).
\end{equation}
From~\eqref{eq_proof_v0_Green_mult_ratios_integral} and~\eqref{eq_proof_v0_Green_mult_ratios_integral_B},
\begin{equation*}
	v_k(x) \le k (b-a)m((a,b)) v_{k-1}(x).
\end{equation*}
Since $(b-a)m((a,b))\leq \xcnorm{X}{\bg}$, we get the desired result on $v_k/v_{k-1}$.
\end{proof}

\begin{corollary}\label{cor_vk_bound}
Let $k\in \mathbb{N} $. Then, we have the following bound for $v_k(x)=\Esp_x(\tau_{ab}^k \indicB{\tau_b< \tau_a})$,
\begin{equation*}
	\|v_k\|_{\infty} \le k! \xcnorm{X}{\bg}^k.
\end{equation*}	
\end{corollary}

\begin{lemma}\label{lem_exp_tau_1}
Let $X$ be a diffusion process with state space $\II $ an interval of $\IR $, $m $ the speed measure of $X$, $a,x,b \in \II $ such that $a<x<b $, $ \bc_x = (a,b)$, $\xcnorm{X}{\bc_x}=m(\bc_x)|b-a| $ and $\lambda > 0 $ such that $\lambda \xcnorm{X}{\bc_x} \in (0,1)$. Then,
\begin{equation}\label{eq_exp_tau_1}
	\Esp_x(e^{\lambda \tauab}) \le  \exp\Big(\frac{\lambda}{1-\lambda \xcnorm{X}{\bc_x}} \Esp_x( \tauab) \Big).
\end{equation}
\end{lemma}

\begin{proof}
Developing the exponential series,
\begin{equation*}
	\Esp_x(e^{\lambda \tauab})
	= 1 + \sum_{N=1}^{\infty} \frac{\lambda^N}{N!} \Esp_x(\tauab^N) = 1 + \Esp_x(\tauab)  \sum_{N=1}^{\infty} \frac{\lambda^N}{N!}  \frac{\Esp_x(\tauab^N)}{\Esp_x(\tauab) }. 
\end{equation*}
From Corollary~\ref{cor_vk_bound}, we can bound the ratio of expected values by $(N-1)!\xcnorm{X}{\bc_x}^{N-1} $. As $\lambda \xcnorm{X}{\bc_x} \in (0,1)$,
\begin{equation*}
	\Esp_x(e^{\lambda \tauab}) \le 1 + \Esp_x(\tauab) \sum_{N=1}^{\infty} \frac{\lambda}{N}(\lambda \xcnorm{X}{\bc_x})^{N-1} \le  1 + \Esp_x(\tauab) \frac{\lambda}{1-\lambda \xcnorm{X}{\bc_x}}.
\end{equation*}
Thus, we only need to apply the inequality $1+x\le e^x $ to get \eqref{eq_exp_tau_1}.
\end{proof}

\begin{lemma}
\label{lem_tau_PM_bound}
Let $t,M>0 $ and $\lambda >0 $ such that  $\lambda \xcnorm{X}{\bg}\in (0,1)  $. Then,
\begin{equation}\label{eq_lem_tau_PM_bound}
	\Prob_x(\tauET{t}{\bg}>M) \le e^{-\lambda M} e^{\lambda \frac{t + \xcnorm{X}{\bg}}{1-\lambda \xcnorm{X}{\bg}}}.
\end{equation}
\end{lemma}

\begin{proof}
From Markov's inequality,
\begin{equation}
	\label{eq_lemma_PM_1}
	\Prob_x(\tauET{t}{\bg}>M) 
	\le e^{-\lambda M} \Esp_x\Big[ e^{\lambda \tauET{t}{\bg}} \Big] 
	= e^{-\lambda M} \Esp_x\Big[ e^{\lambda \sum_{k=1}^{\Kri{t}{\bg}}(\tau_k - \tau_{k-1})} \Big].
\end{equation}
Conditioning on the $\sigma$-algebra $\mathcal{B} $ generated by the trajectory of $X_t $ on the grid $\bg $, \textit{i.e.}, $\mathcal{B} = \sigma\{ X_{\tau_k}; k\in \mathbb{N}_0\} $, and as $\Kri{t}{\bg}$ is $\mathcal{B}$-measurable,
\begin{equation}
	\label{eq_lemma_PM_2}
	\Esp_x\Big[ e^{\lambda \sum_{k=1}^{\Kri{t}{\bg}}(\tau_k - \tau_{k-1})} \Big]
	= \Esp_x\Big[ \prod_{k=1}^{\Kri{t}{\bg}} \Esp\Big[ e^{\lambda(\tau_k -\tau_{k-1})}\Big| \mathcal{B}\Big] \Big].
\end{equation}
From the definition of $\xcnorm{X}{\bg} $, $\lambda  \xcnorm{X}{\bc} \le \lambda \xcnorm{X}{\bg} \in (0,1)$ for each cell $\bc $ of the grid $\bg $. 
Thus, applying Lemma \ref{lem_exp_tau_1} on each term of the product in \eqref{eq_lemma_PM_2},
\begin{equation}
	\label{eq_lem_tau_PM_bound_pre}
	\Prob_x(\tauET{t}{\bg}>M) 
	\le e^{-\lambda M} \Esp_x\Big[  \exp \Big( \frac{\lambda}{1-\lambda \xcnorm{X}{\bg}} \sum_{k=1}^{\Kri{t}{\bg}}  \Esp(\tau_k - \tau_{k-1} | \mathcal{B} ) \Big) \Big]
\end{equation}
From \eqref{eq_random_index_definition},
\begin{equation*}
	t< \sum_{k=1}^{\Kri{t}{\bg}}  \Esp(\tau_k - \tau_{k-1} | \mathcal{B} ) \le 
	t + \Esp(\tau_{\Kri{t}{\bg} +1} - \tau_{\Kri{t}{\bg}} | \mathcal{B} )
\end{equation*}
Thus, from \eqref{eq_vk_ratio_ninf_bound}, 
\begin{equation}\label{eq_bound_sum_embed_times_proof}
	\sum_{k=1}^{\Kri{t}{\bg}}  \Esp(\tau_k - \tau_{k-1} | \mathcal{B} )  \le t + \|v_1 /v_0 \|_{\infty} \le t + \xcnorm{X}{\bg},
\end{equation} 
From \eqref{eq_lem_tau_PM_bound_pre} and \eqref{eq_bound_sum_embed_times_proof}, we get \eqref{eq_lem_tau_PM_bound}.
\end{proof}

\subsection {Convergence of the embedding times}

For this section, let $X $ be a diffusion process with state-space $\II $ an interval of $\IR $, $\bg $ a covering grid of $\II $, $\tau_{k}^{\bg} $ the embedding times of $X$ in $\bg $ as defined in \eqref{eq_tauk_def} and $\Kri{t}{\bg} $ as defined in \eqref{eq_random_index_definition}.

\begin{lemma}\label{lem_var_tau}
For any $T>0 $,
\begin{equation*}
	\sum_{k=1}^{\Kri{T}{\bg} }\Var\big( \tau_{k}^{\bg} - \tau_{k-1}^{\bg} \big| X_{\tau_{k-1}^{\bg}} , X_{\tau_{k}^{\bg}} \big) \le  2 \xcnorm{X}{\bg} \big(  T + \xcnorm{X}{\bg} \big),
\end{equation*}
where $|.|_{X}$ is defined in \eqref{eq_grid_metrics_def}.
\end{lemma}

\begin{proof}
For all $x \in \bg $, let $\bc_x $ be the cell of $\bg $ containing  $x$. Then,
\begin{multline*}
	\sum_{k=1}^{\Kri{T}{\bg} } \Var\big( \tau_{k}^{\bg} - \tau_{k-1}^{\bg} \big| X_{\tau_{k-1}^{\bg}} , X_{\tau_{k}^{\bg}} \big)
	\\  
	=    
	\sum_{k=1}^{\Kri{T}{\bg}  }  
	\CBesp{(\tau_{k}^{\bg} - \tau_{k-1}^{\bg})^2 }{ X_{\tau_{k-1}^{\bg}} , X_{\tau_{k}^{\bg}} }{}
	-\Big(\CBesp{\tau_{k}^{\bg} - \tau_{k-1}^{\bg} }{ X_{\tau_{k-1}^{\bg}} , X_{\tau_{k}^{\bg}} }{} \Big)^2
	\\	
	\le 
	\sup_{k\le \Kri{t}{\bg} } \Bigg\{ \frac{
		\CBesp{(\tau_{k}^{\bg} - \tau_{k-1}^{\bg})^2 }{ X_{\tau_{k-1}^{\bg}} , X_{\tau_{k}^{\bg}} }{}
	}
	{
		\CBesp{\tau_{k}^{\bg} - \tau_{k-1}^{\bg} }{ X_{\tau_{k-1}^{\bg}} , X_{\tau_{k}^{\bg}} }{}
	}\Bigg\}
	\sum_{k=1}^{\Kri{T}{\bg} }  
	\CBesp{(\tau_{k}^{\bg} - \tau_{k-1}^{\bg}) }{ X_{\tau_{k-1}^{\bg}} , X_{\tau_{k}^{\bg}} }{}.
\end{multline*}
So from Lemma~\ref{lem_vk_ratio} and the definition of $\Kri{t}{\bg} $,
\begin{multline*}
	\sum_{k=1}^{\Kri{T}{\bg}  } \Var\big( \tau_{k}^{\bg} - \tau_{k-1}^{\bg} \big| X_{\tau_{k-1}^{\bg}} , X_{\tau_{k}^{\bg}} \big)  \\
	\le 
	\Big\|\frac{v_2}{v_1} \Big\|_{\infty}
	\sum_{k=1}^{\Kri{T}{\bg} }  
	\CBesp{(\tau_{k}^{\bg} - \tau_{k-1}^{\bg}) }{ X_{\tau_{k-1}^{\bg}} , X_{\tau_{k}^{\bg}} }{}
	\\
	\le \Big\|\frac{v_2}{v_1} \Big\|_{\infty} \Big( T + \Big\|\frac{v_1}{v_0} \Big\|_{\infty}  \Big)
	\le 2\xcnorm{X}{\bg} (T + \xcnorm{X}{\bg}),
\end{multline*}
which is the desired inequality.
\end{proof}

\begin{proposition}
\label{prop_tau_cond_mart}
Let $\process{\mathcal{F}_t} $ be the canonical filtration of $X $.
Let also $\mathcal{A}_n = \mathcal{F}_{\tau_{n}^{\bg}} $, $\mathcal{B}= \sigma\big((X_{\tau_{k}^{\bg}})_{k\in\mathbb{N}_0}\big) $ and $\Delta \tau_{k}^{\bg} = \tau_{k}^{\bg} - \tau_{k-1}^{\bg} $. If we define the augmented filtration $\mathcal{G}_n = \mathcal{A}_n \vee \mathcal{B} $, then  the process
\begin{equation*}
	M_n = \sum_{k=1}^{n} \Delta \tau_{k}^{\bg} - \Esp \big[\Delta \tau_{k}^{\bg} \big| X_{\tau_{k-1}^{\bg}} , X_{\tau_{k}^{\bg}} \big],
\end{equation*}
is a $\mathcal{G}_n$-martingale.
\end{proposition}

\begin{proof}
For $m\le n $,
\begin{multline*}
	\Esp \big[M_n \big| \mathcal{G}_m \big] =
	\Esp \big[M_n \big| \mathcal{A}_m,\mathcal{B} \big]\\ = \Esp\Big[\sum_{k=1}^{m} \Delta \tau_{k}^{\bg} - \Esp \big[\Delta \tau_{k}^{\bg} \big| X_{\tau_{k-1}^{\bg}} , X_{\tau_{k}^{\bg}} \big] 
	+ \sum_{k=m+1}^{n} \Delta \tau_{k}^{\bg} - \Esp \big[\Delta \tau_{k}^{\bg} \big| X_{\tau_{k-1}^{\bg}} , X_{\tau_{k}^{\bg}} \big] \Big| \mathcal{A}_m,\mathcal{B} \Big].
\end{multline*}
From Lemma~\ref{lem_sideways_markov_property} and as $\tau_{k}^{\bg} $ is $\mathcal{A}_k $ measurable,
\begin{equation*}
	\Esp \big[M_n \big| \mathcal{G}_m \big] = M_m + \sum_{k=m+1}^{n} \Esp_{X_{\tau_{m}^{\bg}}}\Big[ \Delta \tau_{k}^{\bg} - \Esp \big[\Delta \tau_{k}^{\bg} \big| X_{\tau_{k-1}^{\bg}} , X_{\tau_{k}^{\bg}} \big]  \Big| \mathcal{B} \Big] = M_m.
\end{equation*}
This proves the result.
\end{proof}

\begin{theorem}\label{thm_sup_tau_bound}
For any $T>0$, 
\begin{equation}
	\label{Meq_thm_sup_tau_bound}
	\Esp \Big[\sup_{t\in [0,T]}|\tauET{t}{\bg}      -t|^2 \Big]              
	\le 
	2 \xcnorm{X}{\bg} \Big( 4 \big(T + \xcnorm{X}{\bg}\big)  + 1 \Big),
\end{equation}
where $|.|_{X}$ is defined in \eqref{eq_grid_metrics_def}.
\end{theorem}

\begin{proof}
The convexity inequality $(a+b)^p \le 2^{p-1}(a^p+b^p)$ 
yields for $p=2 $	
\begin{multline}
	\label{eq_tau_decomp}
	\Esp \Big[\sup_{t\in [0,T]}|\tauET{t}{\bg}  -t|^2 \Big] \le 2 \Esp \Big[\sup_{t\in [0,T]}\Big|\tauET{t}{\bg}  -
	\sum_{k=1}^{\Kri{t}{\bg}  }\Esp \big[\Delta \tau_{k}^{\bg} \big| X_{\tau_{k-1}^{\bg}} , X_{\tau_{k}^{\bg}} \big]\Big|^2 \Big]\\
	+ 2 \Esp \Big[\sup_{t\in [0,T]}\Big|\sum_{k=1}^{\Kri{t}{\bg}}\Esp \big[\Delta \tau_{k}^{\bg} \big| X_{\tau_{k-1}^{\bg}} , X_{\tau_{k}^{\bg}} \big] - t\Big|^2 \Big]. 
\end{multline}
From Proposition \ref{prop_tau_cond_mart}, the term $	\sum_{k=1}^{n} \Delta \tau_{k}^{\bg} - \Esp \big[\Delta \tau_{k}^{\bg} \big| X_{\tau_{k-1}^{\bg}} , X_{\tau_{k}^{\bg}} \big]$ is a $\mathcal{G}_n$-martingale, where $\mathcal{G}_n =  \mathcal{F}_{\tau_{n}^{\bg}} \vee \mathcal{B}$ and $\mathcal{B}= \sigma\big((X_{\tau_{k}^{\bg}})_{k\in\mathbb{N}_0}\big)  $.
Thus, from Doob's  $L^p $ inequality,
\begin{equation*}\label{eq_tau_decomp_part1}
	\Esp \Big[ \sup_{t\in[0,T]} |M_{\Kri{t}{\bg}}|^2 \Big] = \Esp \Big[ \sup_{k\le \Kri{T}{\bg}} |M_k|^2 \Big] \le 2 \Esp \big[ |M_{\Kri{T}{\bg}}|^2 \big]. 
\end{equation*}
By conditioning on $\mathcal{B} $, from Lemma \ref{lem_sideways_markov_property},
\begin{multline*}
	\label{eq_tau_decomp_part1_2}
	\Esp \Big[\sup_{t\in [0,T]}\Big|\sum_{k=1}^{\Kri{t}{\bg}} \Delta \tau_{k}^{\bg}- \Esp \big[\Delta \tau_{k}^{\bg} \big| X_{\tau_{k-1}^{\bg}} , X_{\tau_{k}^{\bg}} \big]\Big|^2 \Big]
	\\
	\le 2	\Esp \Big[\Big(\sum_{k=1}^{\Kri{T}{\bg}} \Delta \tau_{k}^{\bg}- \Esp \big[\Delta \tau_{k}^{\bg} \big| X_{\tau_{k-1}^{\bg}} , X_{\tau_{k}^{\bg}} \big]\Big)^2 \Big] 
	\\
	= 2 \Esp \Big[\sum_{k=1}^{\Kri{T}{\bg}} \Var \big(\Delta \tau_{k}^{\bg} \big| X_{\tau_{k-1}^{\bg}} , X_{\tau_{k}^{\bg}}  \big)  \Big],
\end{multline*}
which from Lemma~\ref{lem_var_tau} is bounded by $ 4 \xcnorm{X}{\bg} \big(  T + \xcnorm{X}{\bg} \big) $.
For the second term on the right hand side of~\eqref{eq_tau_decomp}, from~\eqref{eq_random_index_definition} and since $\Kri{t}{\bg} \ge 1 $ for any $t>0 $,
\begin{multline*}
	\sum_{k=1}^{\Kri{t}{\bg}} \Esp \big[\tau_{k}^{\bg} - \tau_{k-1}^{\bg} \big| X_{\tau_{k-1}^{\bg}} , X_{\tau_{k}^{\bg}} \big] - t 
	\le \Esp \big[\tau_{\Kri{t}{\bg}}^{\bg} - \tau_{\Kri{t}{\bg}-1}^{\bg} \big| X_{\tau^{\bg}_{\Kri{t}{\bg} - 1}} , X_{\tauET{t}{\bg}} \big]
	\\
	\le \Big\|\frac{v_1}{v_0} \Big\|_{\infty} \le \xcnorm{X}{\bg}.
\end{multline*}
So having bounded both additive parts of the right hand side of~\eqref{eq_tau_decomp}, we get~\eqref{Meq_thm_sup_tau_bound}.
\end{proof}

\section{Convergence rate of the Markov chain}\label{sec:convergence}

\subsection{Moment bounds}

In this section we prove that Theorem~3.1 of \cite{AnkKruUru2} holds also for reflected processes and for a sharper constant.
This result, combined with the bound \eqref{Meq_thm_sup_tau_bound}
is used to prove the convergence of the approximation process in Section \ref{ssec:main_proof}.

\begin{lemma}\label{lem_moments_diff_bound}
Let $X $ be a diffusion process on natural scale with state-space $\IR$ and a speed measure $m_X $ that satisfies Condition~\eqref{eq_condition1}. Then, for all $p\geq 2$, there exist two constants $C$, $C'>0$ such that
\begin{equation}
	\label{Meq_lem_moments_diff_bound}
	\Esp_x \Big[\sup_{t\in [0,T]} |X_t - x|^p \Big] \le  C'\,[1+|x|^{p}]e^{C T},
\end{equation}
where $k_1$ is a constant such that Condition~\eqref{eq_condition1} is satisfied, $C\le 8p (p-1)/k_1$  and $C'>0 $ is a constant that depends only on $p$.
\end{lemma}

\begin{proof}
Let $Z$ be the diffusion process on natural scale with speed measure,
\begin{equation*}
	m_Z(\d x)= \indicB{|x|<1} \frac{k_1}{2} \vd x + \indicB{|x|\ge 1} \frac{k_1}{2x^2} \vd x.
\end{equation*}	
We note that $\frac{k_1}{2x^2}\leq \frac{k_1}{1+ x^2}$ for all $x\geq 1$. 
The dynamic of $Z $ can be shown to be
\begin{equation*}
	\vd Z_t=
	\begin{cases}
		\dfrac{2}{\sqrt{k_1}}   \vd B_t, &\text{for }|Z_t|< 1, \\
		\dfrac{2}{\sqrt{k_1}}  Z_t \vd B_t,& \text{for }|Z_t|\ge 1,
	\end{cases}
\end{equation*}
where $B $ is a standard Brownian motion. We also assume that~$\Prob_x (Z_0 = x) = 1$.
As~$X$ and~$Z $ are on natural scale, they can be expressed as time-changed Brownian motion
\cite[Theorem~47.1, p.~277]{RogWilV2}, \textit{i.e.}, for every $t\ge 0 $,	$X_t = B_{\gamma_X(t)} $ and $Z_t = W_{ \gamma_Z(t) }$, where $B $ and $ W$ are two standard Brownian motions with $\Prob_x(B_0 = x)=\Prob_x(W_0 = x) = 1$, respective local times $L^{x}(B),L^{x}(W) $ and with $\gamma_X(t), \gamma_Z(t)$  being the respective right-inverses\footnote{The right-inverse of a function $f$ is given by, $$f^{-1}(x)=  \inf\{\zeta\ge 0: f(\zeta)>x \}.$$}  of
\begin{equation*}\label{eq_proof_A_inv_time_changes}
	A_X(t)=\frac{1}{2}\int_{\II} L^{x}_{t}(B)m_X(\d x)
	\qquad \text{and} \qquad 
	A_Z(t)=\frac{1}{2}\int_{\II} L^{x}_{t}(W)m_Z(\d x).
\end{equation*} 
Using the same underlying Brownian motion in these definitions, we have $A_Z(t)\le A_X(t) $, and hence $\gamma_X(t)\le \gamma_Z(t) $. Thus,
\begin{multline}\label{eq_proof_X_bounds_intermsZ_1}
	\Esp_x \Big[\sup_{t\in [0,T]} |X_t - x|^p \Big]= \Esp_x \Big[\sup_{t\in [0, \gamma_X(T)]} |B_t - x|^p \Big] \\
	\le \Esp_x \Big[\sup_{t\in [0, \gamma_Z(T)]} |B_t - x|^p \Big]  = \Esp_x \Big[\sup_{t\in [0,T]} |Z_t - x|^p \Big].
\end{multline}
Thus, from Doob's $L^p $ and convexity inequalities,
\begin{equation}
	\label{eq1_lem_moments_diff_bound}
	\Esp_x \Big[\sup_{t\in [0,T]} |X_t - x|^p \Big] \le \frac{2^{p-1}p}{p-1} \Big( \Esp_x\big[|Z_T|^{p} \big] + |x|^{p} \Big).
\end{equation}
Using Itô's formula, followed by a classical localization argument, Fatou's Lemma along with the standard dominated convergence theorem, one obtains that for all $q> 2/3$,
\begin{multline*}
	\Esp_x\Big[(|Z_t|-1)^{3q}\indic{|Z_t|\geq 1}+1\Big]  
	\\ \leq 
	(|x|-1)^{3q}\indic{|x|\geq 1} +1+\int_0^t \frac{3q}{2}(3q-1)
	\Esp_x\Big[(|Z_s|-1)^{3q-2}\indic{|Z_s|\geq 1}\frac{4Z_s^2}{k_1}\Big]\vd s.
\end{multline*}
Using the inequality $|x-1|^{3q-2}x^{2}\leq  4 |x-1|^{3q}+4$ for all $x\geq 1$,
\begin{multline*}	
	\Esp_x\Bigsqbraces{(|Z_t|-1)^{3q}\indic{|Z_t|\geq 1}+1}\\ 
	\leq (|x|-1)^{3q}\indic{|x|\geq 1}+1+\int_0^t \frac{24q}{k_1}(3q-1) \Esp_x\Bigsqbraces{(|Z_s|-1)^{3q}\indic{|Z_s|\geq 1}+1}\vd s.
\end{multline*}
Using Gronwall's Lemma, we deduce that, for all $t\geq 0$,
\begin{equation*}	
	\Esp_x\Bigsqbraces{(|Z_t|-1)^{3q}\indic{|Z_t|\geq 1}+1}
	\leq [(|x|-1)^{3q}\indic{|x|\geq 1}+1]e^{24q (3q-1) t/k_1}.
\end{equation*}
Let $C'_q>0$ be a constant such that $ (|x|-1)^{3q}\indic{|x|\geq 1}+1 \geq C'_q |x|^{3q}$ for all $x\in \mathbb{R}$. Thus,
\begin{equation*}
	\label{eq2_lem_moments_diff_bound}
	\Esp_x\Big[|Z_t|^{3q}\Big]\leq C'_q\,[(|x|-1)^{3q}\indic{|x|\geq 1}+1]e^{ 24q (3q-1) t/k_1},\quad\forall t\geq 0.
\end{equation*}
Hence, for all\footnote{The extension to $p=2$ is straightforward.} $p> 2$, we have
\begin{equation}
	\label{eq2_lem_moments_diff_bound_B}	
	\Esp_x\Big[|Z_t|^{p}\Big]\leq C'_{p/3}\,[1+|x|^{p}]e^{ 8p (p-1) t/k_1},\quad\forall t\geq 0.
\end{equation}
Replacing \eqref{eq2_lem_moments_diff_bound_B} in \eqref{eq1_lem_moments_diff_bound}, we get the bound \eqref{Meq_lem_moments_diff_bound} for
$C'= \frac{p}{p-1}  2^{p-1} (C'_{p/3} + 1)$ and $C= 8p (p-1)/k_1$.

\end{proof}

\begin{proposition}
Let $X $ be a diffusion process with state-space $\II $ an interval of $\IR $, on natural scale and with a speed measure $m_X $ that satisfies Condition~\eqref{eq_condition1}. Then, for each $T>0 $ and $\gamma\in (0,\frac{1}{2}) $, there exists a constant $C>0$ such that
\begin{equation}
	\label{prop_Kol_cents}
	\biglpnorm{  \sup_{s\ne t \le T} \frac{|X_t - X_s|}{|t-s|^{\gamma}}  }{p}{x} \le C (1 + |x|).
\end{equation}
The result also holds in the presence of a reflecting boundary $\zeta \in \II $.
\end{proposition}

\begin{proof}
For the non-reflecting case, the proof works using the same arguments as in the proof of Theorem~3.1 in~\cite{AnkKruUru2}.

For the reflecting case:
Let $X$ be a diffusion process on natural scale, with speed measure $m$  satisfying~\eqref{condition8} for a constant $k_1>0 $ and a reflecting boundary at $\zeta \in \II$.
We observe that 
\begin{equation}
	\label{eq_proof_reflected_as_abs_plus_X}
	X = |X^{\circ} - \zeta| + \zeta 
\end{equation}
in law, where $X^{\circ} $ is the non-reflecting diffusion on natural scale with speed measure 
\begin{equation*}
	\label{eq_proof_symmetrized_refl_process}
	m_{X^{\circ}}(\d x) = \indicB{x\ge \zeta} m(\d x)
	+ \indicB{x < \zeta}m(2\zeta - \vd x)
\end{equation*}
which also satisfies \eqref{condition8} for the same constant $k_1 $.
From \eqref{eq_proof_reflected_as_abs_plus_X} and the triangle inequality,
\begin{equation}\label{eq_proof_reflecting_main_A}
	\frac{|X_t - X_s|}{|t-s|^{\gamma}} 
	=  \frac{\bigabsbraces{|X^{\circ}_t-\zeta| - |X^{\circ}_s- \zeta|}}{|t-s|^{\gamma}} \le \frac{|X^{\circ}_t - X^{\circ}_s|}{|t-s|^{\gamma}}.  
\end{equation}
The diffusion process $X^{\circ} $ is non-reflecting and its speed measure $m_{X^{\circ}} $ satisfies \eqref{condition8}.
Thus,~\eqref{prop_Kol_cents} holds for $X^{\circ} $ for a constant $C>0 $ and from \eqref{eq_proof_reflecting_main_A}, 
\begin{equation*}
	\label{prop_Kol_cents_refl}
	\biglpnorm{  \sup_{s\ne t \le T} \frac{|X_t - X_s|}{|t-s|^{\gamma}}  }{p}{x}
	\le \biglpnorm{  \sup_{s\ne t \le T} \frac{|X^{\circ}_t - X^{\circ}_s|}{|t-s|^{\gamma}}  }{p}{x} \le C (1 + |x|).
\end{equation*}
\end{proof}

\subsection{Proof of the convergence rate for a process on natural scale}\label{ssec:main_proof}

\begin{proof}[Proof of Theorem \ref{thm_main}]

From Proposition \ref{prop_embedding_2}, 
\begin{multline*}
	\mathcal{W}_p\Big[\Law\big ((\widetilde{X}^{\bg}_{t})_{t\in [0,T]} \big) ,\Law\big ((X_{t})_{t\in [0,T]} \big)\Big]\\
	=         \inf
	\Bigcubraces{ \Big\| 
		\|\zeta - \xi \|
		\Big\|_{L^p}; \zeta \sim \Law\big ((\widetilde{X}^{\bg}_{t})_{t\in [0,T]} \big), \xi \sim \Law\big ((X_{t})_{t\in [0,T]} \big)}\\
	\le       \biglpnorm{  \sup_{t\in [0,T]} |X_{\tauET{t}{\bg}}-X_t  | }{p}{x}.
\end{multline*}
Let $M>T $. From Minkowski inequality,
\begin{multline}
	\label{eq1_thm_main}
	\biglpnorm{\sup_{t\in [0,T]}  |X_{\tauET{t}{\bg}}-X_t  | }{p}{x} \le 
	\biglpnorm{\indicB{\tauET{T}{\bg} \le M} \sup_{t\in [0,T]}  |X_{\tauET{t}{\bg}}-X_t  | }{p}{x} 
	\\
	+ \biglpnorm{\indicB{\tauET{T}{\bg}> M}  \sup_{t\in [0,T]}  |X_{\tauET{t}{\bg}}-X_t  | }{p}{x}
\end{multline}
For the first additive term of \eqref{eq1_thm_main}, for any $\gamma>0 $, by multiplying and dividing by $|\tauET{t}{\bg} - t|^{\gamma} $,
\begin{multline*}
	\Big\|\indicB{\tauET{T}{\bg} \le M} \sup_{t\in [0,T]}  |X_{\tauET{t}{\bg} }-X_t  | \Big\|_{L^p(\Prob_x)} 
	\\
	= 
	\biglpnorm{
		\indicB{\tauET{T}{\bg} \le M}  \sup_{t\in [0,T]} \bigg\{\frac{ |X_{\tauET{t}{\bg} }-X_t  |}{|\tauET{t}{\bg} - t|^{\gamma} } |\tauET{t}{\bg}  - t|^{\gamma}  \bigg\}
	}{p}{x}\\
	\le \bigg\| 
	\indicB{\tauET{T}{\bg} \le M} \sup_{t\in [0,T]} \bigg\{\frac{ |X_{\tauET{t}{\bg} }-X_t  |}{|\tauET{t}{\bg}  - t|^{\gamma} }\bigg\} \sup_{t\in [0,T]} |\tauET{t}{\bg} - t|^{\gamma}   \bigg\|_{L^p(\Prob_x)}.
\end{multline*}
As $\tauET{t}{\bg} $ is increasing with respect to $t$,
\begin{multline*}
	\Big\|\indicB{\tauET{T}{\bg}\le M} \sup_{t\in [0,T]}  |X_{\tauET{t}{\bg}}-X_t  | \Big\|_{L^p(\Prob_x)} \\
	\le \bigg\| \sup_{s\ne t \le M} \bigg\{\frac{ |X_s-X_t  |}{|s - t|^{\gamma} }\bigg\} \sup_{t\in [0,T]} |\tauET{t}{\bg} - t|^{\gamma}   \bigg\|_{L^p(\Prob_x)} \\
	=\biggbraces{ 
		\Esp_x\biggsqbraces{ 
			\sup_{s\ne t \le M} \biggcubraces{ \frac{ |X_s-X_t  |}{|s - t|^{\gamma} }} \sup_{t\in [0,T]} |\tauET{t}{\bg}- t|^{\gamma}  
		}^p 
	}^{1/p}.
\end{multline*}
From Hölder's inequality for $q \ge 1$ and $q/(q - 1)$ conjugates exponents 
and \eqref{prop_Kol_cents}, for every $\gamma \in (0,\frac{1}{2}) $,
there exists a constant $C_1=C_1(M, \gamma, p(q-1)/q,x)>0$ such that
\begin{multline*}
	\Big\| \indicB{\tauET{T}{\bg}\le M} \sup_{t\in [0,T]}  |X_{\tauET{t}{\bg}}-X_t  | \Big\|_{L^p(\Prob_x)} \\
	\le	\bigg( \Esp_x \bigg[  \sup_{s\ne t \le M} \frac{ |X_s-X_t  |}{|s- t|^{\gamma} } \bigg]^{p(q-1)/q} \bigg)^{q/p(q-1)} 
	\big(\Esp_x  \big[ \sup_{t\in [0,T]} |\tauET{t}{\bg}- t|^{\gamma}  \big]^{pq} \big)^{1/pq} \\
	\le   C_{1}^{q/p(q-1)} 
	\big(\Esp_x \big[  \sup_{t\in [0,T]} |\tauET{t}{\bg}  - t|^{\gamma p q} \big]  \big)^{1/pq}.
\end{multline*}
Choosing $q = 2/\gamma p$, 
from \eqref{Meq_thm_sup_tau_bound}, for any $\gamma \in (0, \frac{1}{2}\wedge \frac{2}{p} ) $
\begin{equation}
	\label{eq2_thm_main}
	\Big\| \indicB{\tauET{T}{\bg}\le M} \sup_{t\in [0,T]}  |X_{\tauET{t}{\bg}}-X_t  | \Big\|_{L^p(\Prob_x)} 
	\le C_{1} ^{1/p(1-\gamma p)} \big( 2 \xcnorm{X}{\bg} \left( 4T  + \xcnorm{X}{\bg}\right)\big)^{\gamma/2}.
\end{equation}
For the second additive term of \eqref{eq1_thm_main},
\begin{equation}
	\label{eq3_thm_main}
	\Esp_x \Big[   \indicB{\tauET{T}{\bg}> M}  \sup_{t\in [0,T]}  |X_{\tauET{t}{\bg}}-X_t  |^p \Big] = \sum_{m=  M}^{\infty} \Esp_x \Big[\indicB{\tauET{T}{\bg} \in [m,m+1)}  \sup_{t\in [0,T]}  |X_{\tauET{t}{\bg}}-X_t  |^p  \Big].
\end{equation}
For each term of the sum in \eqref{eq3_thm_main}, from Hölder's inequality,
\begin{multline*}
	\Esp_x \Big[\indicB{\tauET{T}{\bg}\in [m,m+1)}  \sup_{t\in [0,T]}  |X_{\tauET{t}{\bg}}-X_t  |^p  \Big]
	\\  
	\le 
	\Big[\Prob_x (\tauET{T}{\bg}>m)\Big]^{1/q'}
	\Big[ \Esp_x \Big[ \indicB{\tauET{T}{\bg}\in [m,m+1)} \sup_{t\in [0,T]}  |X_{\tauET{t}{\bg}}-X_t  |^p  \Big]\Big]^{q'/(q'-1)}
	\\
	\le 
	\Big[\Prob_x (\tauET{T}{\bg}>m)\Big]^{1/q'}
	\Big[ \Esp_x \Big[ \indicB{\tauET{T}{\bg}<m+1}\sup_{t\in [0,T]}  |X_{\tauET{t}{\bg}}-X_t  |^p  \Big]\Big]^{q'/(q'-1)}.
\end{multline*}
As each $m$ in the sum in \eqref{eq3_thm_main} satisfies $m+1>M>T$, from Minkowski's inequality,
\begin{multline*}
	\Esp_x \Big[ \indicB{\tauET{T}{\bg}<m+1}\sup_{t\in [0,T]}  |X_{\tauET{t}{\bg}}-X_t  |^p  \Big]
	\\
	\le \Esp_x \Big[ \indicB{\tauET{T}{\bg}<m+1} 2^{p-1} \sup_{t\in [0,T]}  \Big\{|X_{\tauET{t}{\bg}}-x|^p +|X_t-x  |^p  \Big\} \Big]
	\\
	\le 2^{p}   \Esp_x \Big[ \sup_{t\in [0,m+1]}  |X_t-x|^p \Big], 
\end{multline*}
which from Lemma~\ref{lem_moments_diff_bound} is bounded by
$C_2\,[1+|x|^{p}]e^{C_3 (m+1)} $, where $C_2$ and $C_3$ are positive
constants depending on $p$. This, along with Lemma
\ref{lem_tau_PM_bound} and Hölder's inequality gives us for 
$\lambda>0 $ chosen such that 
$\alpha = \lambda \xcnorm{X}{\bg}<1 $,
\begin{equation*}
	\Esp_x \Big[\indicB{\tauET{T}{\bg}\in [m,m+1)}  \sup_{t\in [0,T]}  |X_{\tauET{t}{\bg}}-X_t  |^p  \Big] 
	\le 
	\Big[C_2\,[1+|x|^{pq'}]e^{C_3 (m+1)} \Big]^{1/q'}
	\Big[e^{-\lambda m} e^{\lambda \frac{T + \xcnorm{X}{\bg}}{1-\lambda \xcnorm{X}{\bg}}}\Big]^{(q'-1)/q'}
\end{equation*}
where $C_2(pq')\le 2pq' (pq'-1)/c$ and $C_2(pq')>0$ a positive constant depending only on $pq' $. Thus, setting $C_4(pq'):=\Big[C_2(pq')\,[1+|x|^{pq'}] \Big]^{1/q'}e^{C_3/q'}>0$, 
\begin{multline*}
	\Esp_x \Big[\indicB{\tauET{T}{\bg}\in [m,m+1)}  \sup_{t\in [0,T]}  |X_{\tauET{t}{\bg}}-X_t  |^p  \Big] 
	\\
	\le C_4 \exp \Bbraces{ \frac{1}{q'} \Big[ C_3 m+ (q'-1)\Big( \lambda \frac{T + \xcnorm{X}{\bg}}{1-\lambda \xcnorm{X}{\bg}}   -\lambda m \Big) \Big]}
	\\
	= C_4 \exp \Bbraces{ \frac{1}{q'} \Big[ C_3 m + (q'-1)\frac{\alpha}{\xcnorm{X}{\bg}}\Big( \frac{T + \xcnorm{X}{\bg}}{1-\alpha}   - m \Big) \Big]}
	\\
	= C_5 \exp \Bbraces{ \frac{1}{q'} \Big[ C_3 m + (q'-1)\frac{\alpha}{\xcnorm{X}{\bg}}\Big( \frac{T}{1-\alpha}   - m \Big) \Big]},
\end{multline*}
where $C_5 := C_4 \exp(\frac{q'-1}{q'}\frac{\alpha}{1-\alpha}) $.
If we choose $q'>1$ such that\footnote{We remark that if this is satisfied for a grid $\bg $, then it is satisfied for all grids $\bg ' $ such that $\xcnorm{X}{\bg'}\le \xcnorm{X}{\bg} $.} $A:= C_3-(q'-1) \alpha/ \xcnorm{X}{\bg}<0$, 
\begin{multline*}
	\sum_{m=  M}^{\infty} \Esp_x \Big[\indicB{\tauET{T}{\bg}\in [m,m+1)}  \sup_{t\in [0,T]}  |X_{\tauET{t}{\bg}}-X_t  |^p  \Big] 
	\\
	\le  C_5 \exp\Bbraces{\frac{(q'-1)\alpha T}{(1-\alpha)\xcnorm{X}{\bg}}}  \sum_{m=  M}^{\infty} (e^{A})^m 
	= C_5 \exp\Bbraces{\frac{(q'-1)\alpha T}{(1-\alpha)\xcnorm{X}{\bg}}} \frac{e^{AM}}{1-e^{A}}
	\\
	= C_5 \exp\Bbraces{\frac{(q'-1)\alpha T}{(1-\alpha)\xcnorm{X}{\bg}}} \frac{\exp\Bbraces{\Big(  C_3-\frac{(q'-1) \alpha}{ \xcnorm{X}{\bg}} \Big)M}}{1-\exp\Bbraces{  C_3-\frac{(q'-1) \alpha}{ \xcnorm{X}{\bg}}}}.
\end{multline*}
If $M$ is chosen such that $M>T/(1-\alpha) $, then
\begin{equation}
	\label{eq10_thm_main}
	\sum_{m=  M}^{\infty} \Esp_x \Big[\indicB{\tauET{T}{\bg}\in [m,m+1)}  \sup_{t\in [0,T]}  |X_{\tauET{t}{\bg}}-X_t  |^p  \Big] \le \frac{C_5 e^{C_3 M}}{1-e^A} \exp\Bbraces{\frac{(q'-1)\alpha}{\xcnorm{X}{\bg}}\Big(\frac{T}{1-\alpha}-M\Big)},
\end{equation}
where the bound is $O(e^{-1/\xcnorm{X}{\bg}}) $, and can be rewritten as $C^{(1)}e^{-C^{(2)}/\xcnorm{X}{\bg}} $. From \eqref{eq1_thm_main},  \eqref{eq2_thm_main} and~\eqref{eq10_thm_main},
\begin{equation*}
	\label{eq_bound_xing}
	\biglpnorm{\sup_{t\in [0,T]}  |X_{\tauET{t}{\bg}}-X_t| }{p}{x} \le  C_{1} ^{1/p(1-\gamma p)} \big( 2 \xcnorm{X}{\bg} \Big( 4T  + \xcnorm{X}{\bg}\Big)\big)^{\gamma/2}  + C^{(1)}e^{-C^{(2)}/\xcnorm{X}{\bg}}.
\end{equation*}
As $\xcnorm{X}{\bg}^{\gamma} $ and $e^{-1/\xcnorm{X}{\bg}} $ are both $\Ord(\xcnorm{X}{\bg}^{\gamma/2} ) $, there exists a constant $C >0$, such that,
\begin{equation*}
	\biglpnorm{\sup_{t\in [0,T]}  |X_{\tauET{t}{\bg}}-X_t| }{p}{x} \le C \xcnorm{X}{\bg}^{\gamma/2},
\end{equation*}
which is \eqref{eq_thm_main}.
\end{proof}

\section{Computations for the classical SDE case and beyond}\label{sec:setup}

In order to implement Algorithm~\ref{algorithm_3} for simulating paths of a
diffusion $X$, one needs two things.  First, a grid $\bg $ adapted to the scale
function and speed measure of $X$.  Second, good approximations of the
transition probabilities and conditional transition times of $X$ over $\bg $.
The first point was covered in Section \ref{ssec:tuning}.

In this section, we show how to compute the quantities \eqref{algo_quantities}
in the pure SDE case.  The extension to SDE solutions with point-wise
singularities is straightforward.  This allows us via Algorithm
\ref{algorithm_3} to simulate all such processes.

Let $(\mu,\sigma) $ two real-valued functions that satisfy Condition \ref{cond_sde_condition} and $X$ be the diffusion that solves
\begin{equation*}
\label{eq_sde_sol_2}
\d X_t = \mu (X_t)\vd t + \sigma(X_t)\vd B_t,
\end{equation*}
where $B$ is a standard Brownian motion.
Let also $\II $ be the state-space of $X$.
A straightforward computation using Itô's formula gives us the infinitesimal generator of $X $,
\begin{equation}
\label{eq_SDE_ig}
(\Lop ,\dom (\Lop) )=
\begin{cases}
	\Lop f(x) = \mu (x) f'(x) + \frac{1}{2}\sigma^2(x)f''(x), \qquad \forall f \in \dom (\Lop),\\
	\dom (\Lop)=\bigbraces{f\in C_b(\II): \Lop f \in C(\II)}. \\
\end{cases}
\end{equation}
In particular, if $\mu$ and $\sigma $ are continuous, then, $\dom(\Lop) = C^{2}(\II) $.
From Proposition~\ref{prop_v0_DP} and~\eqref{eq_SDE_ig} we can infer that, 
\begin{equation}
\label{eq_dm}
s'(x)= e^{-\int_{a}^{y}\frac{2\mu (\zeta)}{\sigma^2 (\zeta)}\vd \zeta}
\qquad\text{ and }\qquad
m(\d x) = \frac{1}{s'(x)}\frac{2}{\sigma^2(x)}\vd x.
\end{equation}
Thus, from Proposition \ref{prop_v0_DP},
\begin{equation}
\label{eq_SDE_v0}
v_0(x) = \frac{\int_{a}^{x} e^{-\int_{a}^{y}\frac{2\mu (\zeta)}{\sigma^2 (\zeta)}\vd\zeta}\vd y}
{\int_{a}^{b} e^{-\int_{a}^{y}\frac{2\mu (\zeta)}{\sigma^2 (\zeta)}\vd\zeta}\vd y}.
\end{equation}
From \eqref{meq_pro_vk_equiv} and \eqref{eq_dm},
\begin{equation}
\label{eq_integral_form_2}
\begin{split}
	v_1(x) &=\int_{a}^{b}G_{a,b}(x,\zeta) v_0(\zeta)  \frac{1}{s'(x)}\frac{2}{\sigma^2(x)}\vd\zeta, \\
	\overline{v}_1(x) &=\int_{a}^{b}G_{a,b}(x,\zeta) \big(1-v_0(\zeta) \big)  \frac{1}{s'(x)}\frac{2}{\sigma^2(x)}\vd\zeta. 
\end{split}
\end{equation}

\noindent
The scale functions and the speed measures are defined up to a multiplicative 
constant: for $\alpha\in\mathbb{R}$ and $\lambda>0$,
the pairs $(s,m)$ and $(\alpha+\lambda s,\lambda^{-1}m)$ are associated to the same diffusion.
In particular, as we are only concerned with points $x\in[a,b]$, we could use $v_0$
for the scale function. The speed measure shall be adapted accordingly.
For $x,\zeta \in[a,b]$, the Green function in \eqref{eq_Green_fun} takes the simpler form,
\begin{equation*}
G_{a,b}(x,\zeta)=v_0(x\wedge \zeta)\big(1-v_0(x\vee \zeta)\big). 
\end{equation*}
Thus, expressions \eqref{eq_integral_form_2} become
\begin{equation}
\label{eq_integral_form_3}
\begin{split}
	v_1(x) 
	&=\int_{a}^{b} v_0(x\wedge \zeta)\big(1-v_0(x\vee \zeta)\big)\frac{v_0(\zeta)}{v_0'(\zeta)}\frac{2}{\sigma^2(\zeta)}\d\zeta, \\
	\overline{v}_1(x) 
	&=\int_{a}^{b} v_0(x\wedge \zeta)\big(1-v_0(x\vee \zeta)\big)\frac{1-v_0(\zeta)}{v_0'(\zeta)}\frac{2}{\sigma^2(\zeta)}\d\zeta. 
\end{split}
\end{equation}
\noindent
The quantities \eqref{eq_SDE_v0} and \eqref{eq_integral_form_2} or \eqref{eq_integral_form_3}
can be computed analytically as we do for the Ornstein-Uhlenbeck process in Section~\ref{example_OU}
or approximated numerically as for the Cox-Ingersoll-Ross process in Section~\ref{example_OU}.

\section{Numerical experiments}
\label{sec:example}

\subsection{Approximation in distribution}\label{ssec_lawapprox}

%

\subsubsection{Standard Brownian motion}\label{example_BM}

The standard Brownian motion can be defined as (see \cite[p. 119]{BorSal}) the diffusion process with scale function and speed measure
\begin{equation*}
s(x)= x \quad \text{and} \quad
m(\d x)=
2 \vd x.
\end{equation*}
Let $v_k(x):=\Esp_x[\tauab^k \indicB{\tau_b < \tau_a}  ] $ and $\overline{v}_k(x):=\Esp_x[\tauab^k \indicB{\tau_a < \tau_b}  ] $ for all $k\in \mathbb{N}_0 $, where $\bc_x = (a,b) $. Then, from the definition of the scale function,
\begin{equation*}
v_0(x)
=\frac{x-a}{b-a}.
\end{equation*}
From Proposition \ref{pro_vk_equiv},
\begin{equation*}
v^{BM}_1(x)= 
(x-a)(b-x)\Big( \frac{2}{3}\frac{(x-a)^2}{(b-a)^2} + \frac{b-x}{b-a} - \frac{2}{3}\frac{(b-x)^2}{(b-a)^2}   \Big) ,
\end{equation*}
\begin{equation*}
\overline{v}^{BM}_1(x)=
(x-a)(b-x)\Big( \frac{2}{3}\frac{(b-x)^2}{(b-a)^2} + \frac{x-a}{b-a} -\frac{2}{3}\frac{(x-a)^2}{(b-a)^2}   \Big),
\end{equation*}
where by $BM$ we mean these quantities are associated with the Brownian motion.
Thus, from~\eqref{algo_quantities}, we have all the necessary quantities we need to implement the algorithm.

\subsubsection{Sticky Brownian motion}\label{example_stqBM}

The Brownian motion with a sticky point at $0$ of stickiness $\rho >0$, 
already presented in Example~\ref{ex:sticky}, can be
defined \cite[p.~123]{BorSal} as the diffusion process with scale
function and speed measure 
\begin{equation*}
s(x)= x \quad \text{and} \quad
m(\d x)=
2 \vd x + \rho \delta_0 (\d x).
\end{equation*}

\noindent
From the definition of the scale function,
\begin{equation*}
v_0(x) =\frac{x-a}{b-a}.
\end{equation*}
From Proposition \ref{pro_vk_equiv} we may deduce the following expressions for the conditional exit times of $\bc_x = (a,b) $,
\begin{equation*}
\begin{split}
	v_1(x)&= 
	v^{BM}_1(x) + \rho \indic{0 \in (a,b)}  G_{(a,b)}(x,0) v_0(0),\\
	\overline{v}_1(x)&=
	\overline{v}^{BM}_1(x) +  \rho \indic{0 \in (a,b)}  G_{(a,b)}(x,0) \big( 1 - v_0(0) \big),
\end{split}
\end{equation*}
where  $v^{BM}_1(x) $ and $\overline{v}^{BM}_1(x)  $ are the analogous quantities for the standard Brownian motion.

\begin{figure}[h!]
\centering
\begin{subfigure}{0.45\textwidth}
	\scalebox{0.57}{\includegraphics{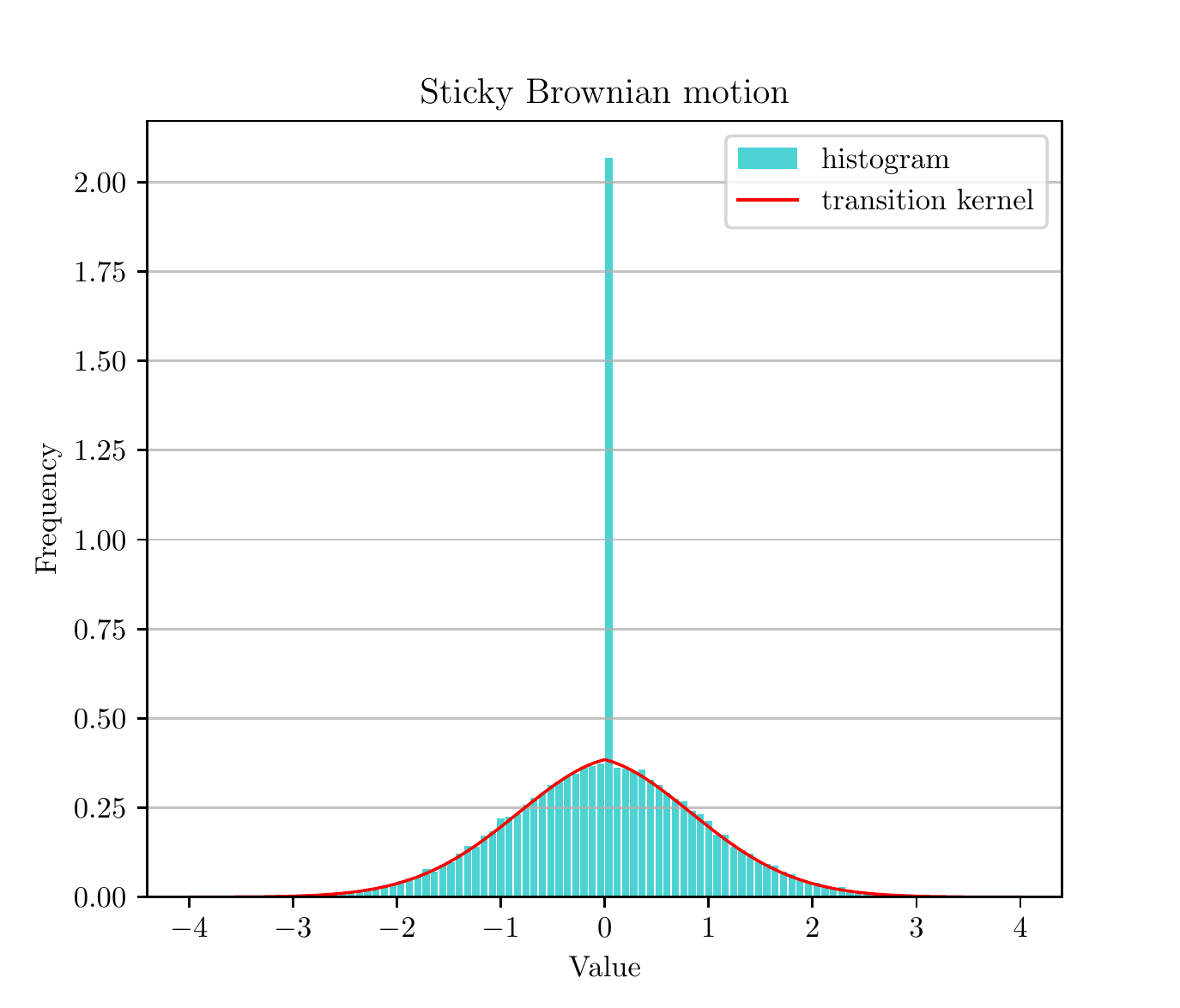}}
	\caption{\label{fig:sticky_pdf}
	}
\end{subfigure}
\qquad
\begin{subfigure}{0.45\textwidth}
	\scalebox{0.57}{\includegraphics{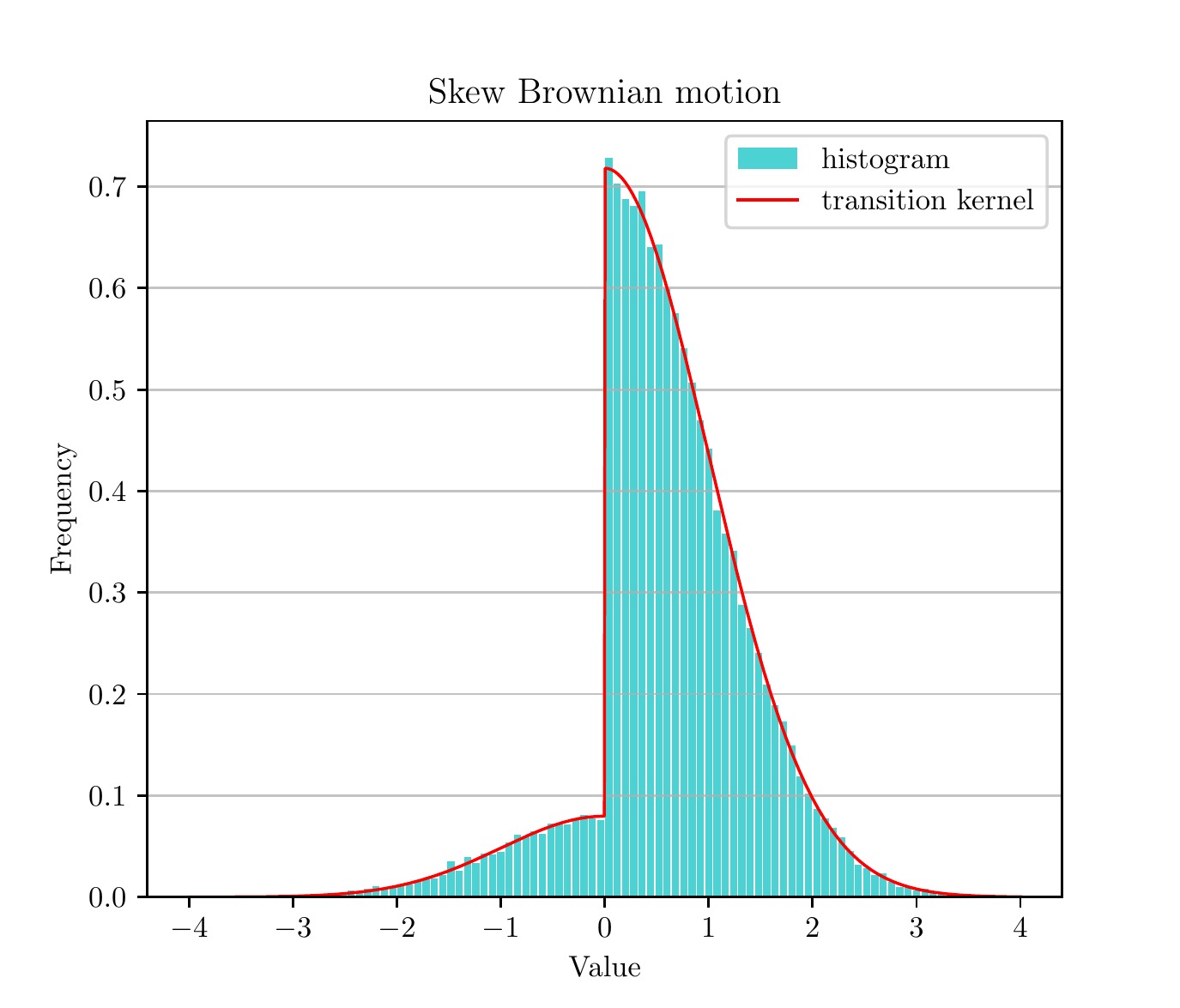}}
	\caption{\label{fig:skew_pdf}
	}
\end{subfigure}
\caption{
	\label{fig:grid_sticky_skew_pdf} 
	(a):
	Histogram of simulated values at $T=1 $ of a sticky Brownian motion of parameter $\rho = 0.7 $ with initial value $x_0 = 0 $ using Algorithm \ref{algorithm_3} with the tuned grid \eqref{eq_tuned_grid} of size-criteria $h=0.01 $.
	\newline
	(b):
	Histogram of simulated values at $T=1 $ of a skew Brownian motion of parameter $\beta = 0.9 $ with initial value $x_0 = 0 $ using Algorithm \ref{algorithm_3} with a uniform grid of step-size $h=0.01 $.
} 
\end{figure}

In Figure \ref{fig:grid_sticky_skew_pdf}-(a) are simulation results of the sticky Brownian 
motion using algorithm \ref{algorithm_3} and its probability transition kernel (see \cite[p. 124]{BorSal}).

\subsubsection{Skew Brownian motion}\label{example_skewBM}

The skew Brownian motion at $0$ with parameter $\beta \in (0,1) $  can be defined (see \cite[p. 126]{BorSal}) as the diffusion process with scale function and speed measure 
\begin{equation*}
\scale (x)=
\begin{cases}
	\dfrac{x}{\beta}, &x\ge 0,\\
	\dfrac{x}{1-\beta},&x\le 0\\
\end{cases}\quad \text{and} \quad
\speed (\d x)=
\begin{cases}
	2\beta \vd x, &x>0,\\
	2(1-\beta) \vd x,&x<0.
\end{cases}
\end{equation*}

From the definition of the scale function,
\begin{equation*}
v_0(x) =\frac{s(x)-s(a)}{s(b)-s(a)}.
\end{equation*}
As the $\beta $ and $(1-\beta)$ terms between the speed measure and the scale function compensate themselves in the expressions giving $v_k $ in  Proposition \ref{pro_vk_equiv},
\begin{equation*}
v_1(x)= v^{BM}_1(x)
\text{ and }
\overline{v}_1(x)= \overline{v}^{BM}_1(x).
\end{equation*}

In Figure \ref{fig:grid_sticky_skew_pdf}-(b) are simulation results of the skew Brownian motion using Algorithm~\ref{algorithm_3} and a plot of its probability transition kernel \cite[p. 126]{BorSal}.

\subsubsection{Ornstein-Uhlenbeck process}\label{example_OU}

Let $X $ be the Ornstein-Uhlenbeck with mean reversion force $\theta >0 $, long-term mean $\mu $ and diffusion parameter $\sigma >0 $ define through its path-wise description
\begin{equation*}
\vd X_t = \theta (\mu - X_t)\vd t + \sigma \vd B_t, \quad X_0\in \mathbb R,
\end{equation*} 
where $B $ is a standard Brownian motion. We will see that $v_0 $ can be expressed in terms of the Gaussian imaginary error function as $s(x)=\erfi \big( \sqrt{\frac{\theta}{\sigma ^2}} (\mu - x)  \big)$, 
\begin{equation*}
v_0(x) =  \frac{\int_{a}^{x} e^{-\int_{a}^{y}\frac{2\mu (\zeta)}{\sigma^2 (\zeta)}d\zeta}\vd y}
{\int_{a}^{b} e^{-\int_{a}^{y}\frac{2\mu (\zeta)}{\sigma^2 (\zeta)}d\zeta}\vd y}
= \frac{\erfi \big( \sqrt{\frac{\theta}{\sigma ^2}} (\mu - x)  \big)  -  \erfi \big( \sqrt{\frac{\theta}{\sigma ^2}} (\mu - a)  \big) }
{ \erfi \big( \sqrt{\frac{\theta}{\sigma ^2}} (\mu - b)  \big) - \erfi \big( \sqrt{\frac{\theta}{\sigma ^2}} (\mu - a)  \big)},
\end{equation*}	
where $\erfi (x) 	= \frac{2}{\pi} \int_{0}^{x}e^{t^2}\vd t = \frac{2}{\pi} e^{x^2}D_+(x)$, with $ D_+(x)$ being the Dawson function\footnote{We use the latter representation in our numerical results.}. 
Thus, for an Ornstein-Uhlenbeck process \eqref{eq_integral_form_2} becomes
\begin{equation*}
\label{eq_integral_form_ou}
\begin{split}
	v_1(x) 
	&=\int_{a}^{b} v_0(x \wedge \zeta) \big(1- v_0(x \vee \zeta)\big) v_0(\zeta) c\exp\Big(\frac{\theta(\zeta - \mu)^2}{\sigma^2}\Big) \frac{2}{\sigma^2}\vd x, \\
	\bar{v}_1(x) 
	&=\int_{a}^{b} v_0(x \wedge \zeta) \big(1- v_0(x \vee \zeta)\big)  \big(1 - v_0(\zeta)\big) c \exp\Big(\frac{\theta(\zeta - \mu)^2}{\sigma^2}\Big)\frac{2}{\sigma^2}\vd x, \\
\end{split}
\end{equation*}
where $c =  \erfi \big( \sqrt{\frac{\theta}{\sigma ^2}} (\mu - b)  \big) - \erfi \big( \sqrt{\frac{\theta}{\sigma ^2}} (\mu - a)  \big)$.

\subsubsection{Cox-Ingersoll-Ross process}
\label{example_CIR}

\begin{figure}[h!]
\centering
\begin{subfigure}{0.45\textwidth}
	\scalebox{0.57}{\includegraphics{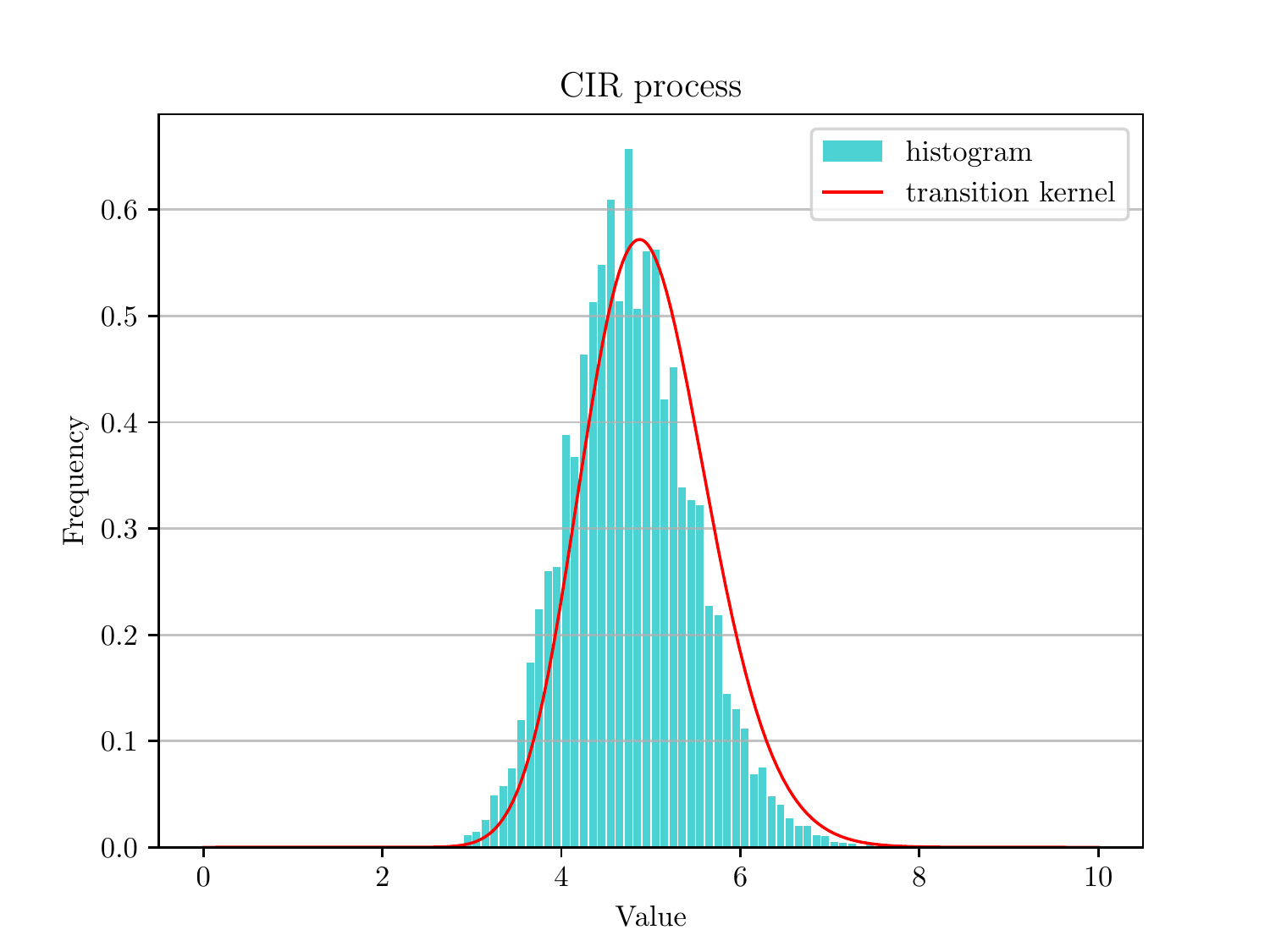}}
	\caption{}
	\label{fig:cir3_uni}
\end{subfigure}
\qquad
\begin{subfigure}{0.45\textwidth}
	\scalebox{0.57}{\includegraphics{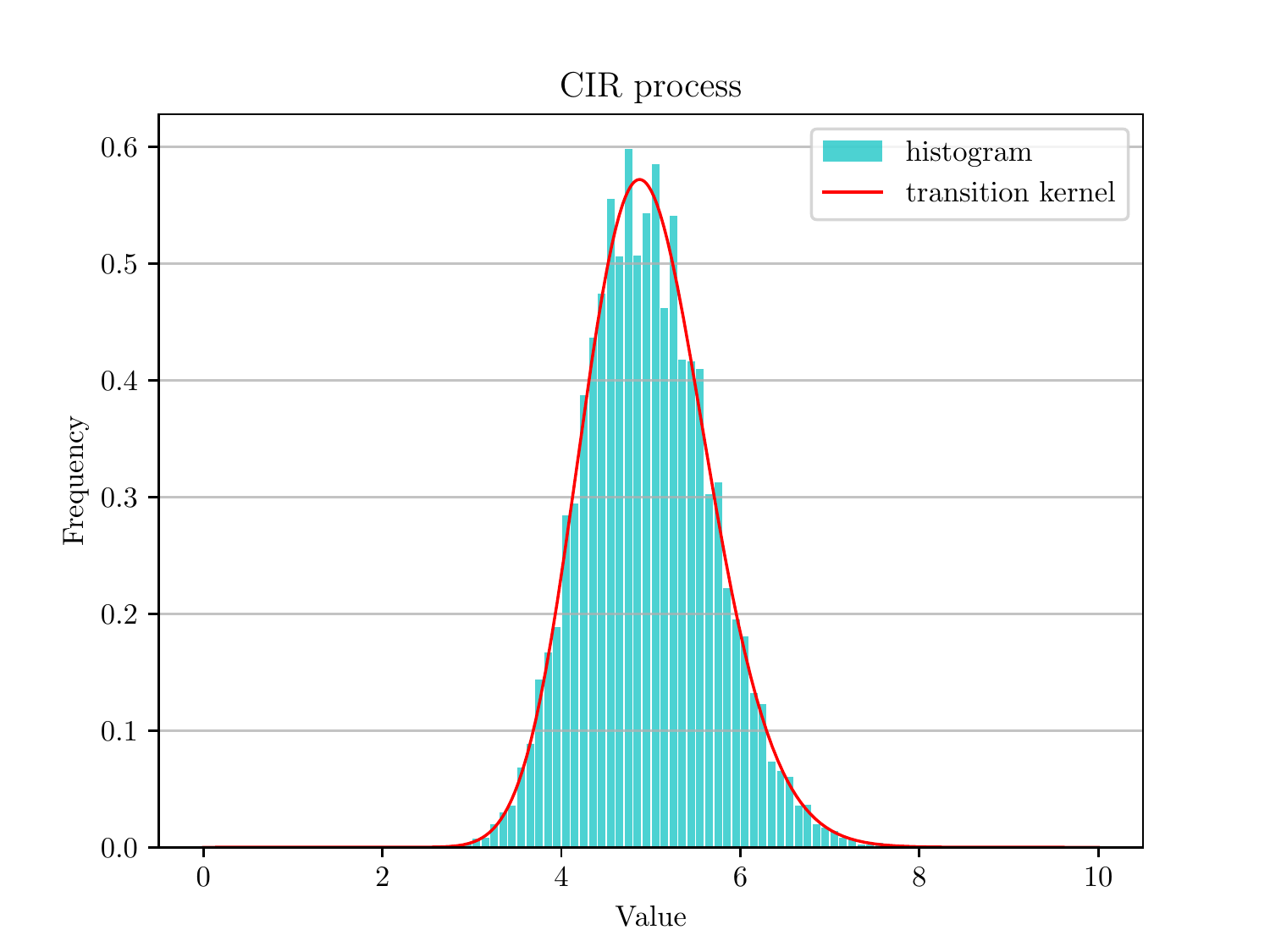}}
	\caption{}
	\label{fig:cir3_tn}
\end{subfigure}
\caption{Histogram of simulated values at $T=1 $ of a CIR process  of $(\theta,\mu,\sigma)=(5,5,1) $ with initial value $x_0 = 1 $ using Algorithm \ref{algorithm_3} with:
	\newline  (a): a uniform grid of step-size $h=0.01 $ and $ (250,200)$-step Riemann approximation of $(v_0,v_1)$ (simulation time: 44.5 sec).
	\newline  (b): a tuned grid of size-criteria $h=0.01 $ computed solving numerically  \eqref{eq_grid_tuning_criteria_SDE}-\eqref{eq_grid_tuning_reiteration} with Newton's method and $ (250,200)$-step Riemann approximation of $(v_0,v_1)$ (simulation time: 46.5 sec).} 
\label{fig:grid_tuning_cir3} 
\end{figure}

\begin{figure}[h!]
\centering
\begin{subfigure}{0.45\textwidth}
	\scalebox{0.57}{\includegraphics{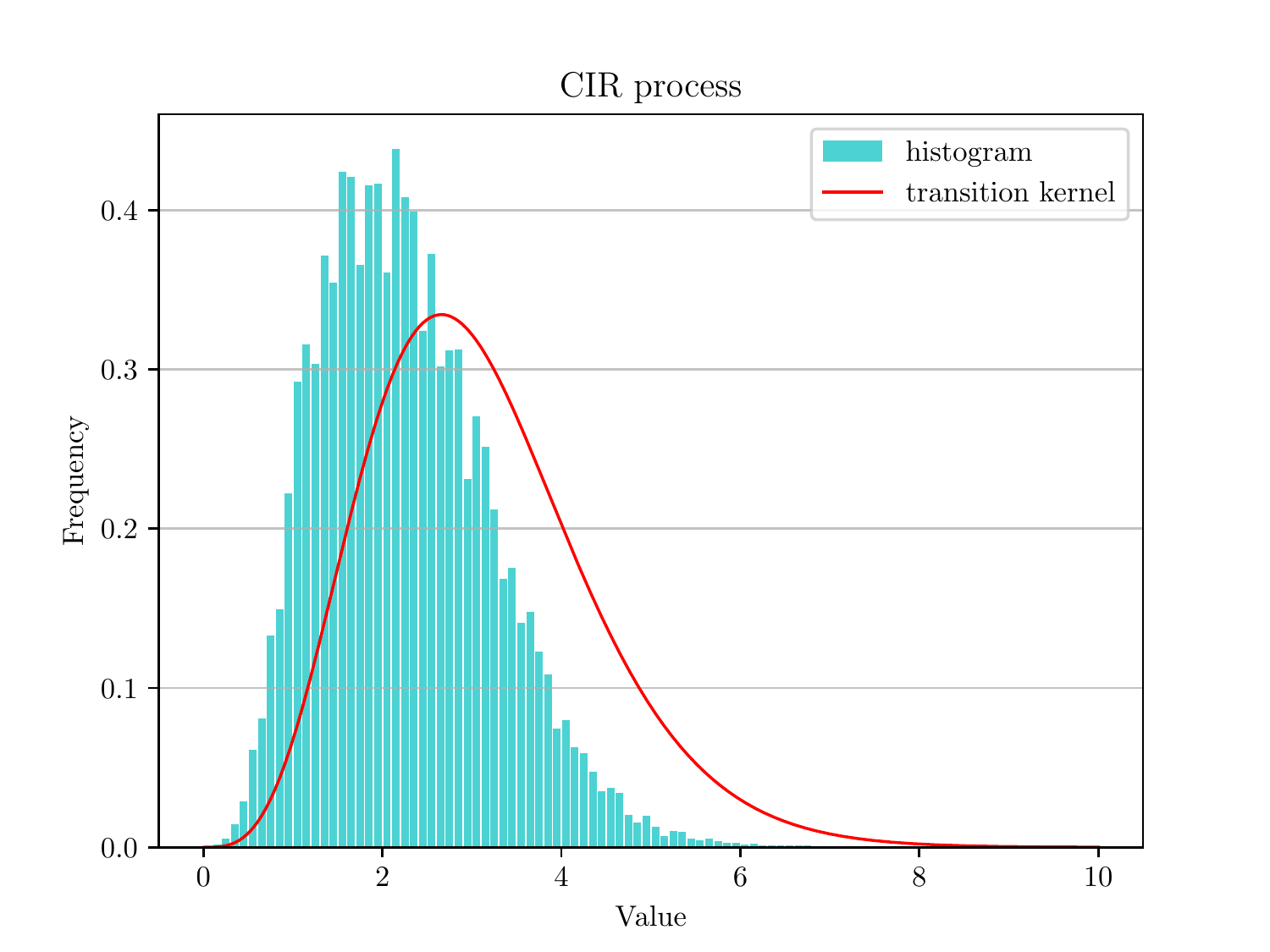}}
	\caption{\label{fig:first}}
\end{subfigure}
\qquad
\begin{subfigure}{0.45\textwidth}
	\scalebox{0.57}{\includegraphics{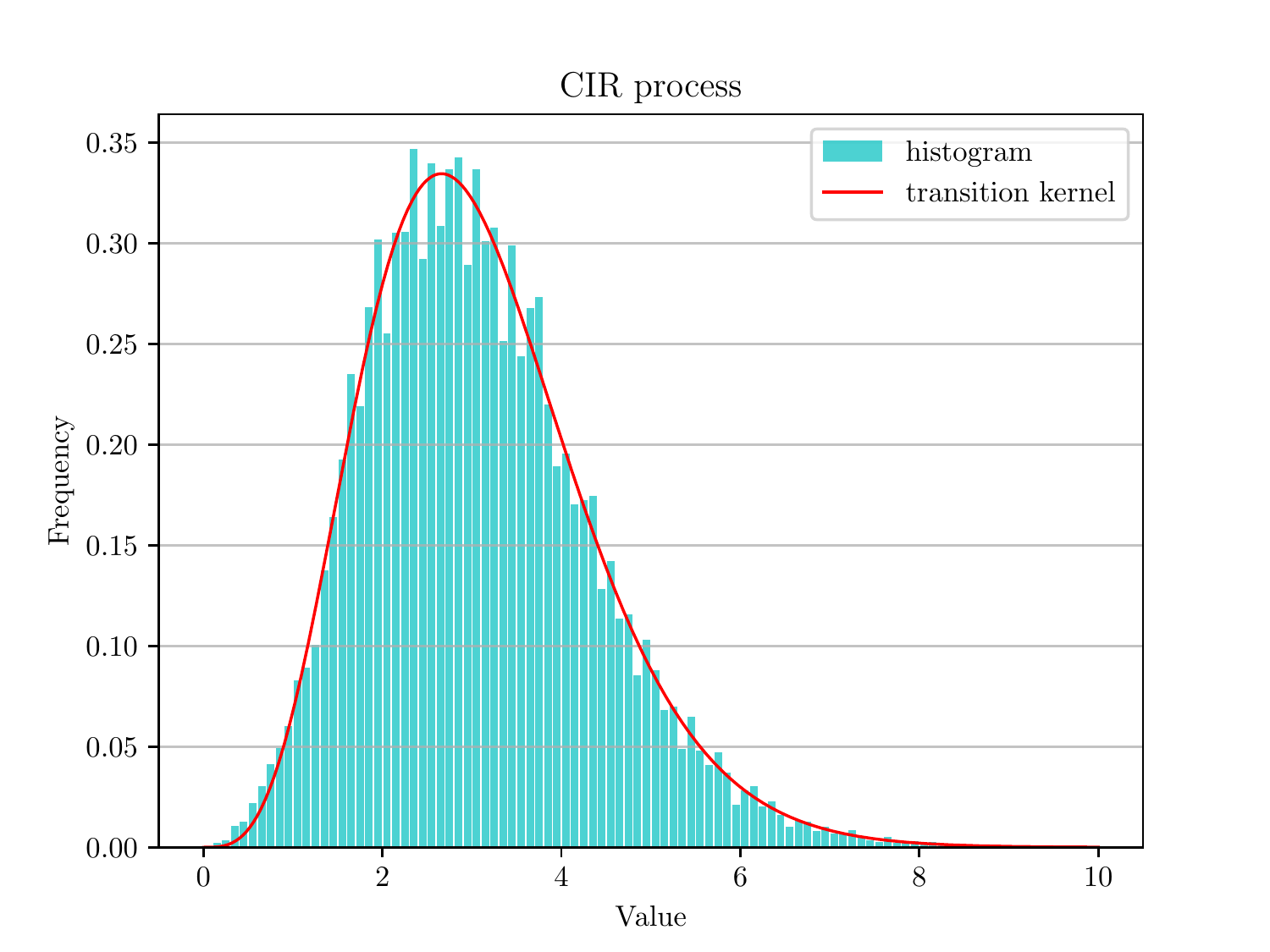}}
	\caption{\label{fig:second}}
\end{subfigure}
\hfill
\begin{subfigure}{0.45\textwidth}
	\scalebox{0.57}{\includegraphics{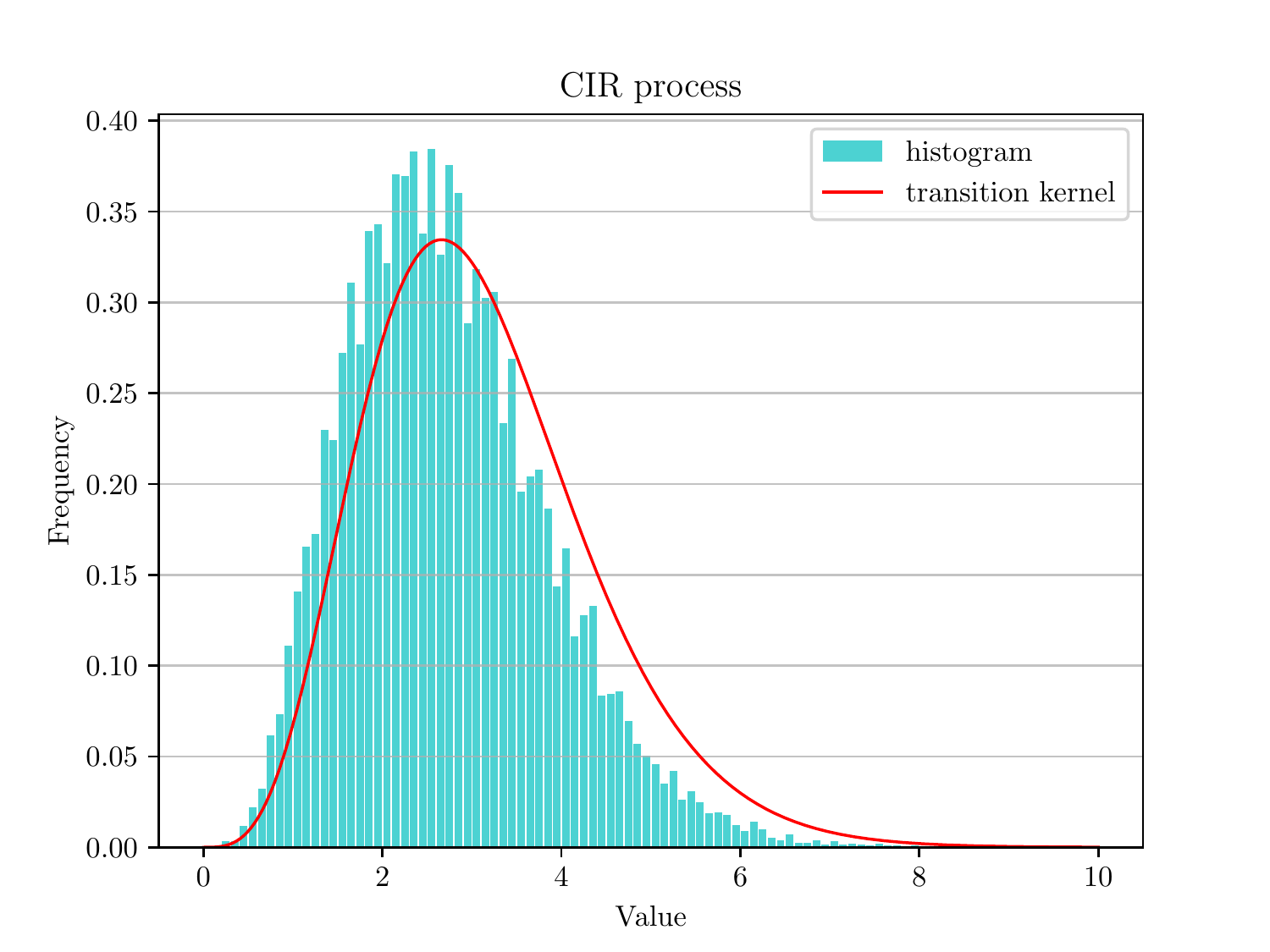}}
	\caption{\label{fig:third}}
\end{subfigure}
\qquad
\begin{subfigure}{0.45\textwidth}
	\scalebox{0.57}{\includegraphics{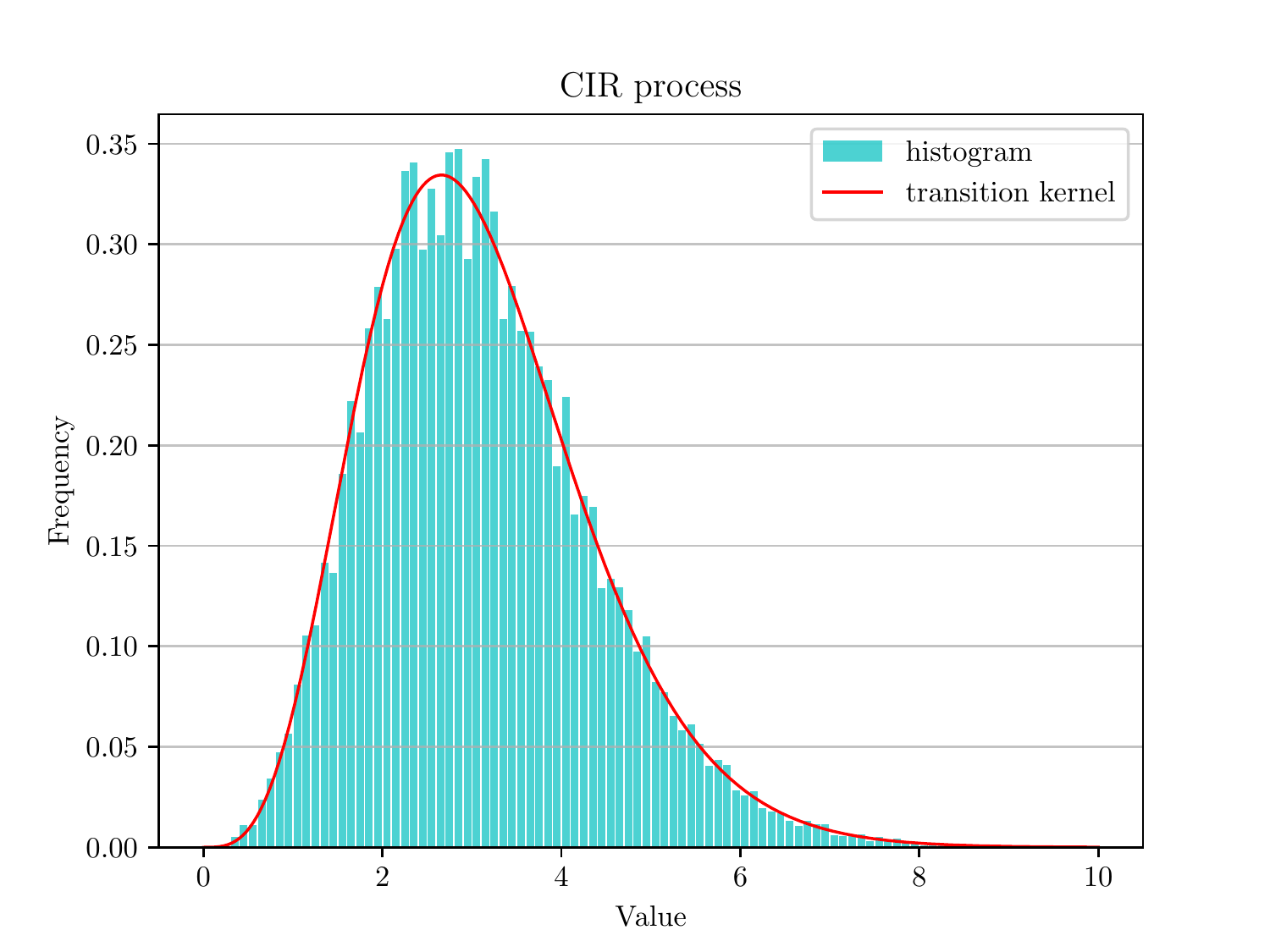}}
	\caption{\label{fig:fourth}}
\end{subfigure}		
\caption{
	\label{fig:grid_tuning}
	Histogram of simulated values at $T=1 $ of a CIR process  of $(\theta,\mu,\sigma)=(1,2,1) $ with initial value $x_0 = 5 $ using Algorithm \ref{algorithm_3} with:
	\newline  (a): a uniform grid of step-size $h=0.015 $ and $ (100,100)$-step Riemann approximation of $(v_0,v_1)$ (simulation time: 10.8 seconds).
	\newline  (b): a tuned grid of size-criteria $h=0.015 $ computed solving numerically  \eqref{eq_grid_tuning_criteria_SDE}-\eqref{eq_grid_tuning_reiteration} with Newton's method and $ (100,100)$-step Riemann approximation of $(v_0,v_1)$ (simulation time: 12.8 seconds). 
	\newline  (c): Same as Figure \ref{fig:first} but with a $(250,200)$-step Riemann approximation of $(v_0,v_1)$ (simulation time: 11.5 seconds).
	\newline  (d): Same as Figure \ref{fig:second} but with a $(250,200)$-step Riemann approximation of $(v_0,v_1)$ (simulation time: 12.2 seconds). }

\end{figure}

The Cox-Ingersoll-Ross process or CIR process \cite{cir}, introduced first by W.~Feller~\cite{feller}, is the diffusion that solves the SDE
\begin{equation*}
\label{eq_cir_process}
\d X_t = \theta (\mu - X_t)\vd t + \sigma \sqrt{X_t}\vd B_t,\  X_0>0,
\end{equation*}
where $B_t$ is a standard Brownian motion. 
The parameter $\theta >0$ expresses its mean reversion speed, $\mu $ its long
term speed and $\sigma >0 $ is its diffusivity parameter.  This equation has a
diffusion coefficients that degenerates at $0$. It however remains non-negative
given $X_0> 0$ and almost surely never hit $0$ when $2\theta\mu>\sigma^2$. A
large body of work have been devoted to the simulation of the CIR and related
process, see e.g.~\cite{alfonsi}. 

From \eqref{eq_SDE_v0} and \eqref{eq_integral_form_3}, 
\begin{equation*}
\label{eq_v0_CIR}
v_0(x) = \frac{\int_{a}^{x} y^{-\frac{2\theta \mu }{\sigma^2}}e^{\frac{2\theta}{\sigma^2}y} \vd y }
{\int_{a}^{b} y^{-\frac{2\theta \mu }{\sigma^2}}e^{\frac{2\theta}{\sigma^2}y} \vd y}
\end{equation*}
and
\begin{equation*}
\label{eq_v1_CIR}
\begin{split}
	v_1(x) 
	&=\int_{a}^{b}v_0(x\wedge \zeta)\big(1-v_0(x\vee \zeta)\big) \frac{v_0(\zeta)}{v_0'(\zeta)}\frac{2}{\sigma^2 \zeta}\vd\zeta, 
	\\
	\overline{v}_1(x) 
	&=\int_{a}^{b}v_0(x\wedge \zeta)\big(1-v_0(x\vee \zeta)\big) \frac{1-v_0(\zeta)}{v_0'(\zeta)}\frac{2}{\sigma^2 \zeta }\vd\zeta. 
\end{split}
\end{equation*}
As the scale function yields no satisfactory closed formula, we perform a numerical approximation of both $v_0$ and the couple $(v_1 , \overline{v}_1)$.

These functions may be computed numerically. One may choose a suitable grid
when the process is close to $0$, leading to improved convergence results. It
is also noteworthy that the approximating process does not cross $0$, a problem
which arise when using Euler type schemes. 

In Figures \ref{fig:grid_tuning_cir3} and \ref{fig:grid_tuning} we plot the
transition kernels and histograms of simulated values of CIR processes at $T=1
$ using Algorithm~\ref{algorithm_3} with uniform and tuned grids and
approximation precisions of the quantities $(v_0,v_1,\bar{v}_1) $.  We also
superpose on the same figures the corresponding transition kernels.

\subsubsection{Skew or reflected Bessel process}\label{example_SBES}

\begin{figure}[h!]
\centering
\scalebox{0.8}{\includegraphics{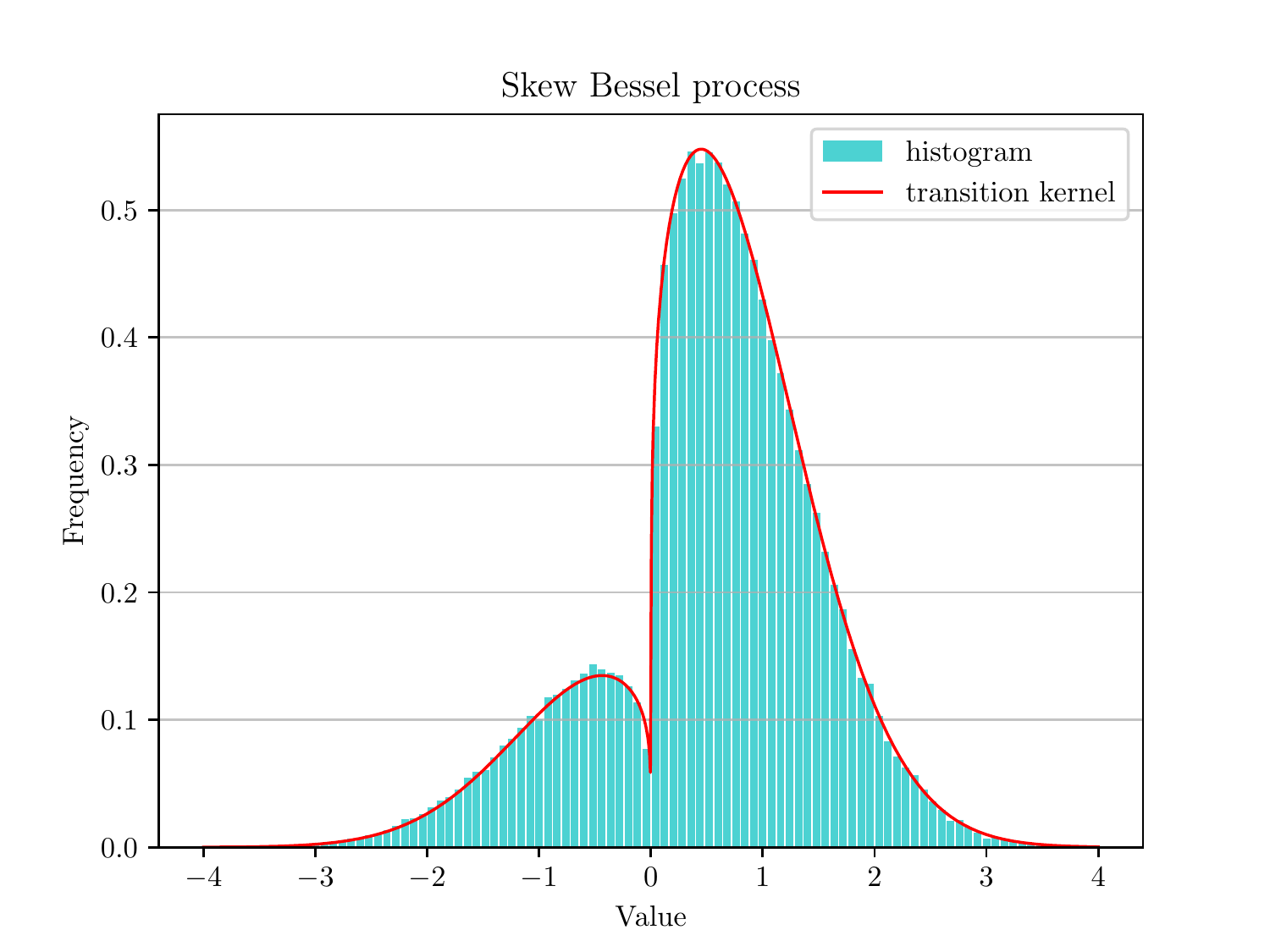}}
\caption{
	\label{fig:SBP}
	Histogram of simulated values at $T=1 $ of the Skew Bessel process of parameters $(\delta,\beta) = (1.2,0.8) $ with initial value $x_0 = 0 $ using Algorithm \protect\ref{algorithm_3} with a uniform grid of step-size $h=0.01 $.
	The quantities \protect\eqref{eq_v0_v1_bessel_integral} were approximated using   100-step Riemann approximations.
}	
\end{figure}

\begin{figure}[h!]
\centering
\begin{subfigure}{0.45\textwidth}
	\scalebox{0.57}{\includegraphics{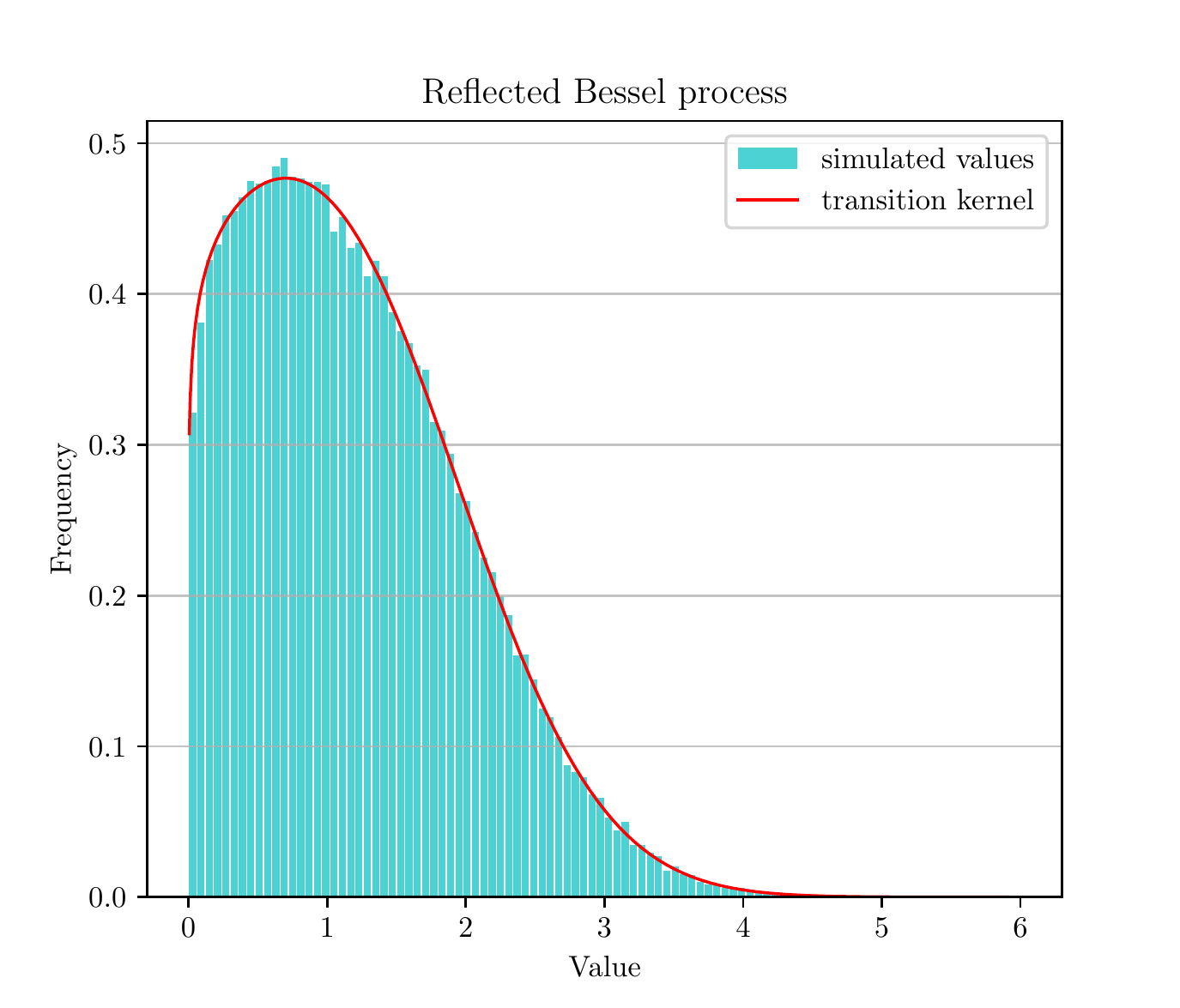}}
	\caption{\label{fig:ssb3_pdf}}
\end{subfigure}
\qquad
\begin{subfigure}{0.45\textwidth}
	\scalebox{0.57}{\includegraphics{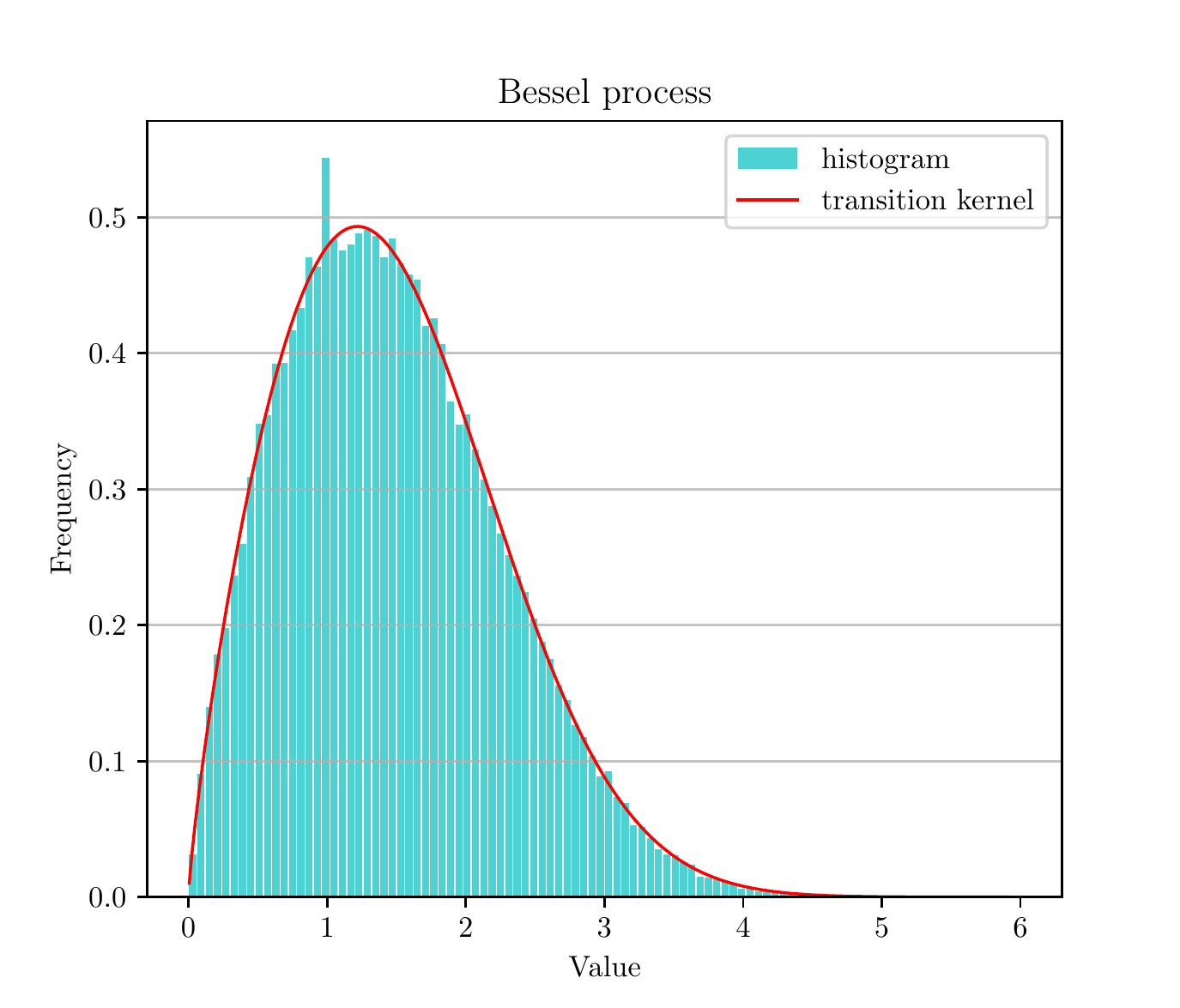}}
	\caption{\label{fig:ssb4_pdf}}
\end{subfigure}
\caption{
	\label{fig:grid_ssb_pdf} 
	(a): histogram of simulated values at $T=1 $ of a Reflected Bessel process of parameter $\delta = 1.1 $ with initial value $x_0 = 5 $ using Algorithm \ref{algorithm_3} with a uniform grid of step-size $h=0.01 $.\protect\newline
	(b): histogram of simulated values at $T=1 $ of a Bessel process of parameter $\delta = 1.8 $ with initial value $x_0 = 5 $ using Algorithm \ref{algorithm_3} with a uniform grid of step-size $h=0.01 $.
	The quantities \eqref{eq_v0_v1_bessel_integral} were approximated using   100-step Riemann approximations.} 
\end{figure}

The Bessel process is defined as the solution of the SDE
\begin{equation*}
\vd X_t = \frac{\delta-1}{2 X_t} \vd t + \vd B_t,
\end{equation*}
where $B$ is a standard Brownian motion.
If $\delta \in (0,2)$, then $0$ is an attainable boundary (see Section \ref{ssec:boundary}).
In this case, we can then define the following behaviors at $0$: reflection or skew/partial reflection. 
In particular, the Skew-Bessel process  of dimension $\delta\in (0,2) $ and skew $\beta \in (0,1) $ is the diffusion process with state-space $\IR $ defined through $s$ and $m$, where (see \cite{AliAyl})
\begin{equation*}
s(x)=  	 
\begin{cases}
	\dfrac{1}{\beta}\dfrac{x^{2-\delta}}{2-\delta}, &x>0,\\
	-\dfrac{1}{1-\beta}\dfrac{|x|^{2-\delta}}{2-\delta}, &x\le 0,\\
\end{cases}
\qquad m(\d x)=			
\begin{cases}
	2\beta x^{\delta - 1}\vd x, &x>0,\\
	2(1-\beta)|x|^{\delta - 1}\vd x, &x\le 0.\\
\end{cases}
\end{equation*}

\noindent
This yield the following expressions for the quantities we need to compute in order to implement the algorithm,
\begin{equation}\label{eq_v0_v1_bessel_integral}
\begin{split}
	v_0(x)
	&= \frac{s(x)-s(a)}{s(b)-s(a)}, \\
	v_1(x) 
	&=\int_{a}^{b}G_{a,b}(x,\zeta) v_0(\zeta)
	2 |\zeta|^{\delta - 1}\vd\zeta, \\
	\overline{v}_1(x) 
	&=\int_{a}^{b}G_{a,b}(x,\zeta) 
	\big(1-v_0(\zeta)\big) 
	2 |\zeta|^{\delta - 1}\vd\zeta, \\
\end{split}
\end{equation}
where $G_{a,b}(x,\zeta)  $ is defined in \eqref{eq_Green_fun}.

In Figures \ref{fig:SBP} and \ref{fig:grid_ssb_pdf} we plot the transition kernels and histograms of simulated values of: a skew and two reflected Bessel processes.
The plotted corresponding transition kernels are computed in \cite{AliAyl}.

\subsection{Local time approximation}\label{ssec_localtime}

The following example illustrates the flexibility of Space-Time Markov Chain Approximation (STMCA) generated via Algorithm~\ref{algorithm_3}.
One feature of such approximation processes is that they are defined on a given grid.
With a suitable choice of grid it is possible to achieve higher orders of convergence of localized path-sensitive functionals.
An example is the stickiness parameter estimation in \cite{Ana22}.
It is shown there that, if $X$ is a sticky Brownian motion of stickiness $\rho>0 $,  $\alpha \in (0,1/2) $ and $g$ is an integrable function vanishing in the vicinity of $0$, the statistics 
\begin{equation}\label{eq_loct_statistic}
\widehat{L}_{n}[g,\alpha] := \frac{n^{\alpha}}{n} \sum_{i=1}^{[nt]} 			g(n^{\alpha}\hfprocess{X}{}{i-1} )
\end{equation}
and
\begin{equation}\label{eq_stickiness_statistic}
\widehat{\rho}_{n}[g,\alpha] := 2   \frac{\lambda(g)}{
	\widehat{L}_{n}[g,\alpha] } \frac{1}{n} \sum_{i=1}^{[nt]} \indic{\bigcubraces{\hfprocess{X}{}{i-1}=0}}
\end{equation}
converge in probability to $\loct{X}{0}{t} $
and $\rho $ respectively.
This qualifies $\widehat{\rho}_{n}[g,\alpha] $
as a consistent estimator of the stickiness parameter.

It is also observed that, as long as $n$ is high enough:
\begin{itemize}
\item The convergence of \eqref{eq_stickiness_statistic}
seems to hold for any $\alpha \in (0,1) $.
\item  The convergence speed of \eqref{eq_stickiness_statistic}
seems to increase in terms of $\alpha \in (0,1) $.
\item The higher the $\alpha$, the less observations of the trajectories are observed via $g$.
\end{itemize}

Algorithm~\ref{algorithm_3} gives us the flexibility to remediate to the latter by using grids of higher precision around the point of stickiness.  Thus, the statistics \eqref{eq_loct_statistic} and \eqref{eq_stickiness_statistic} are relevant for larger values of $\alpha $ and we achieve higher orders of convergence without significant increase in the numerical complexity.

\begin{figure}[h]
\centering
\includegraphics[scale=0.9]{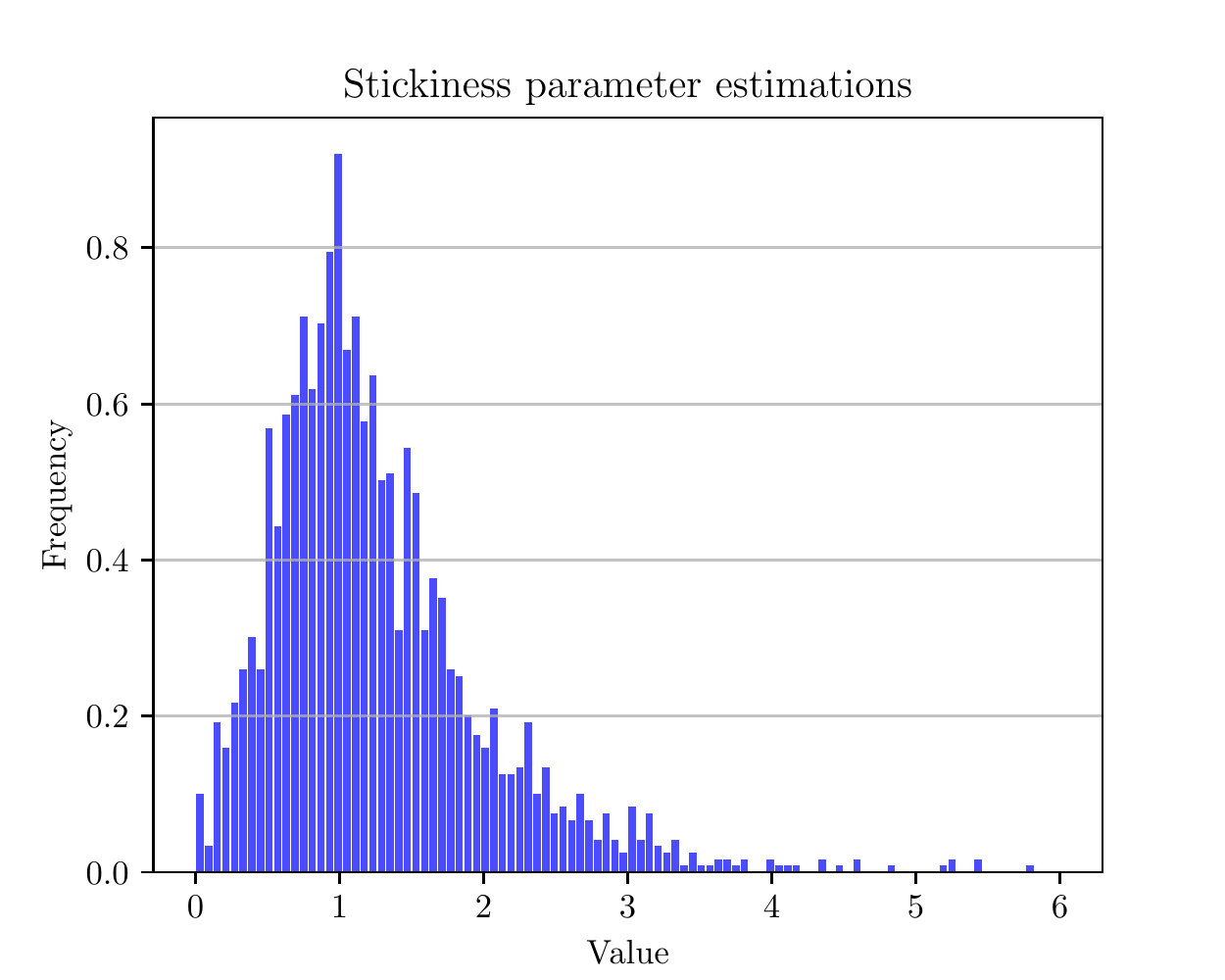}
\caption{\label{fig:se_lowalpha}
	Stickiness parameter estimations histograms using \eqref{eq_stickiness_statistic} with $n = \num{100000} $ and $\alpha = 0.3$ (true value $\rho = 1 $).}
\end{figure}
\begin{figure}[h]
\centering
\includegraphics[scale=0.9]{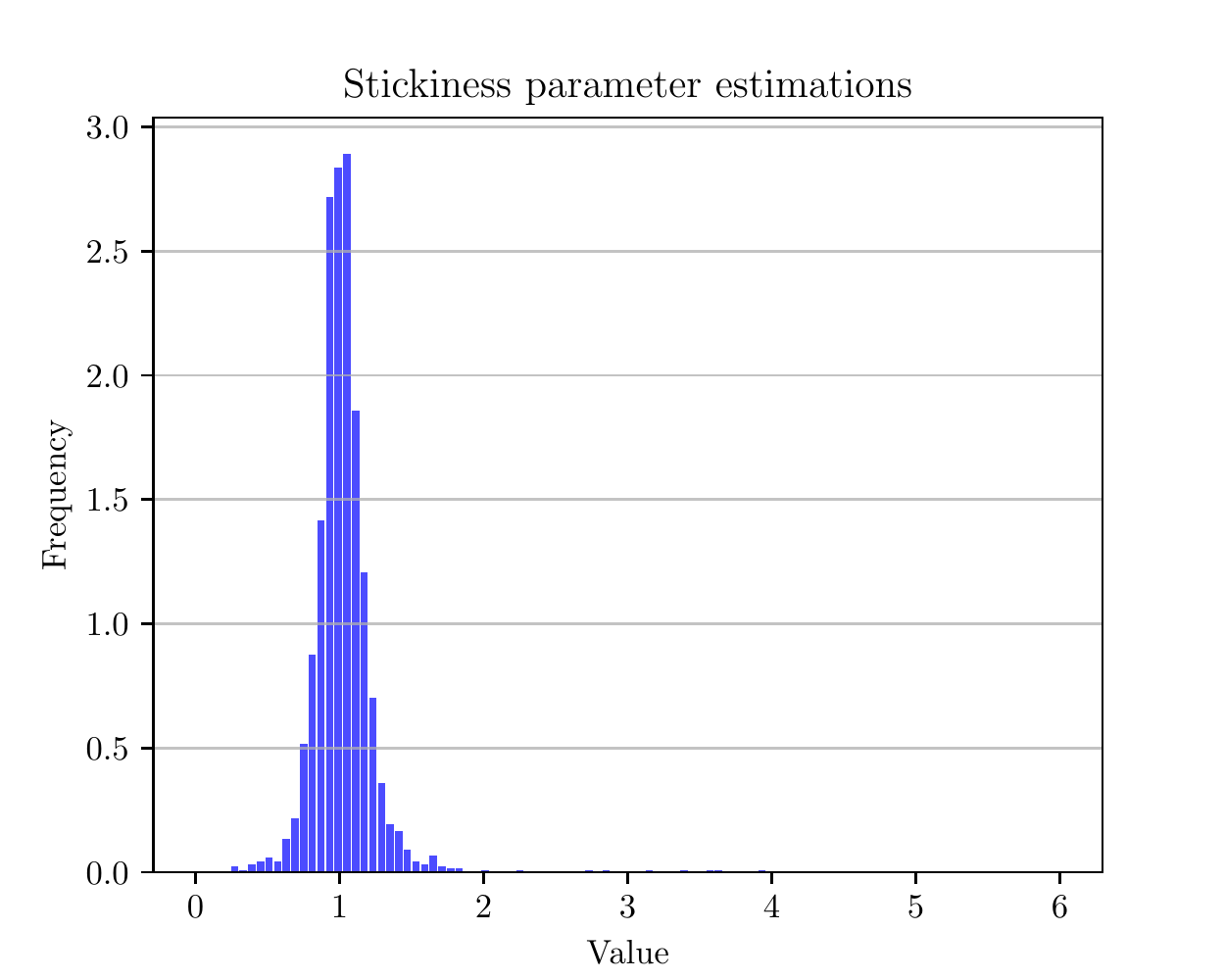}
\caption{
	\label{fig:se_highalpha}
	Stickiness parameter estimations histograms using \eqref{eq_stickiness_statistic} with $n = \num{100000} $ and $\alpha = 0.55$ (true value $\rho = 1 $).}
\end{figure}

Let $\bg_0 $ and $\bg_1 $ be two grids defined by:
\begin{equation}
\bg_0(h) = \{0\} \cup \bigcubraces{\pm (h^2/\rho + k h) ; k \in \IN},\qquad 
\bg_1(h) =  \bigcubraces{\pm x_k(h) ; k \ge 0 },
\end{equation}
where $\{x_j(h)\}_{j\ge 0} $ is defined recursively by
\begin{equation}
x_0 = 0, \qquad
x_j = x_{j-1} + \biggbraces{ \frac{h^2}{\rho} \frac{1}{x_{j-1}+1} + h \Bigbraces{1 - \frac{1}{x_{j-1}+1}}} \indic{x_{j-1} < 1} + h \indic{x_{j-1} \ge 1}.
\end{equation}
We observe that 
\begin{equation}
|\bg_1| = |\bg_0| \qquad  \text{and} \qquad \xcnorm{X}{\bg_1} = \xcnorm{X}{\bg_0}. 
\end{equation}
Thus, from Theorem \ref{thm_main}, the rate of convergence of the STMCA is $ \Ord(h^{1/2})$ for both grids.

\subsubsection*{Simulation results}

\begin{table*}[h!]
	\begin{center}
	\begin{tabular}{@{}llrrrrc@{}}
		\hline
		\multicolumn{1}{c}{$\alpha$} & \multicolumn{1}{c}{$n$} & \multicolumn{1}{c}{$\widehat{\rho}_{\MC}$}
		& \multicolumn{1}{c}{$\widehat{S}_{\MC}^2$} & \multicolumn{1}{c}{$\widehat{\sigma}_{\MC}$}
		& \multicolumn{1}{c}{$\widehat{\acc}$} & \multicolumn{1}{c}{$\mbox{rej}/N_{\MC}$} \\
		\hline
		0.3 & 100000 & 1.334 & 0.740 & 0.860 & 0.125 & 0/2000 \\
		0.4 & 100000 & 1.061 & 0.142 & 0.378 & 0.046 & 0/2000 \\
		0.5 & 100000 & 1.314 & 0.162 & 0.403 & 0.011 & 0/2000 \\
		0.55 & 100000 & -- &  --&  --&  --& 2000/2000 \\
		0.6 & 100000 &  --&  --&  --&  --& 2000/2000 \\
		0.65 & 100000 &  --&  --&  --&  --& 2000/2000 \\
		\hline
	\end{tabular}
	\caption{Stickiness parameter estimations using the grid $\bg_0 $ for $h = 0.01 $. The missing values in the table corresponds to cases where the statistic~\eqref{eq_loct_statistic} is observed to be $0$.  Computation time (single-core): 4 seconds.}
\label{table:non-tuned}
\end{center}
\end{table*}

\begin{table*}[h!]
	\begin{center}
	\begin{tabular}{@{}llrrrrc@{}}
		\hline
		\multicolumn{1}{c}{$\alpha$} & \multicolumn{1}{c}{$n$} & \multicolumn{1}{c}{$\widehat{\rho}_{\MC}$}
		& \multicolumn{1}{c}{$\widehat{S}_{\MC}^2$} & \multicolumn{1}{c}{$\widehat{\sigma}_{\MC}$}
		& \multicolumn{1}{c}{$\widehat{\acc}$} & \multicolumn{1}{c}{$\mbox{rej}/N_{\MC}$} \\
		\hline
		0.3 & 100000 & 1.252 & 0.629 & 0.793 & 0.1272 & 0/2000 \\
		0.4 & 100000 & 1.068 & 0.187 & 0.432 & 0.0437 & 0/2000 \\
		0.5 & 100000 & 1.023 & 0.072 & 0.268 & 0.0140 & 0/2000 \\
		0.55 & 100000 & 1.018 & 0.067 & 0.259 & 0.0079 & 0/2000 \\
		0.6 & 100000 & 1.012 & 0.031 & 0.177 & 0.0044& 0/2000 \\
		0.65 & 100000 & 1.015 & 0.026 & 0.163 & 0.0025 & 2/2000 \\
		\hline
	\end{tabular}
	\caption{Stickiness parameter estimations using the grid $\bg_1$ for $h =0.01 $. Computation time (single-core): 67 seconds.}
\label{table:tuned}
\end{center}
\end{table*}

In this section we present Monte Carlo simulation results.
The integer $N_{\MC}$ will be the Monte Carlo simulation size.
For every $j\le  N_{\MC}$ we simulate the path of an approximation process~$X^{j}_t $.
To assess the quality of each Monte Carlo estimation we use the following metrics:
\begin{itemize}
\item $\bigbraces{\widehat{\rho}_{\MC}, \widehat{S}_{\MC}^2, \widehat{\sigma}_{\MC} }$: Monte Carlo estimation, variance and standard deviation of the stickiness parameter estimations,
\item $\widehat{\acc}$: average number of path-values observed by $g$, i.e
\begin{equation}
	\widehat{\acc} = \frac{1}{N_{\MC}} \sum_{j=1}^{N_{\MC}} \frac{1}{n}\sum_{i=1}^{n}\indic{g(X^{j}_{(i-1)/n}) \ne 0},
\end{equation}
\item rej: percentage of trajectories where the local time estimation equals $0$, \textit{i.e.},
\begin{equation}
	\rej =	
	\# \Bigcubraces{j \le N_{\MC}:
		\frac{n^{\alpha}}{n} \sum_{i=1}^{[nt]} g(n^{\alpha}\hfprocess{X}{j}{i-1} ) = 0
	}.
\end{equation}
\end{itemize}

We use the test function $g(x) = \indic{1<|x|<5}/8 $ which satisfies the conditions of Corollary~1.3 of~\cite{Ana22}.
Within each table (Tables \ref{table:non-tuned}-\ref{table:tuned}) we use the same simulated STMCA trajectories of the sticky Brownian motion of parameter $\rho=1 $.

We observe that the usage of grid $\bg_1 $ yields far superior results than
$\bg_0 $ (see Tables \ref{table:non-tuned}-\ref{table:tuned}).  Using
$\bg_1  $ we have an abundance of simulated path-wise observations close to the
point of stickiness.  The statistic \eqref{eq_stickiness_statistic} remains
thus relevant for large values of $\alpha $ and we can achieve higher orders of
convergence (see Figures~\ref{fig:se_lowalpha}-\ref{fig:se_highalpha}, and
Tables~\ref{table:non-tuned}-\ref{table:tuned}).
It is worth noting that simulation with $h=0.001 $ 


\appendix

\section{Conditioning on the embedded path}

\begin{lemma}\label{lem_sideways_markov_property}
	Let $X$ be a diffusion process with state-space $\II $ an interval of $\IR $
	defined on a family of probability spaces $\bP = (\Omega, \process{\bF_t},\Prob_x)_{x\in \II} $ 
	such that for every $x\in \II $, $\Prob_x (X_0 = x)=1 $.
	Let also $\bg $ be a covering grid of $\II $, $\{\tau_j\}_{j \in \IZ_+} $ the sequence of embedding times of $X$ in $\bg $ defined in \eqref{eq_tauk_def}
	and $\mathcal{B} $  the sigma-algebra defined by
	$\mathcal{B} = \sigma\{ X_{\tau_j}; j\in \mathbb{N}_0\}$.
	Then, for any measurable path-functional $F: C^{0}(\IR_{+}, \II)\mapsto \IR $ and $j\ge 1 $ and $x\in \II $,
	\begin{equation}\label{eq_conditional_equality_inf}
		\Esp_x \bigbraces{ F(\sandwitchprocessusage{X}{j}) \big| \bB} 
		= \Esp_x \bigbraces{ F(\sandwitchprocessusage{X}{j}) \big| 
			X_{\tau_{j-1}},X_{\tau_j} },
	\end{equation}
	where $X^{\tau_{j-1},\tau_{j}} $ is the process defined for every $t\ge 0 $ by $X^{\tau_{j-1},\tau_{j}}_t = X_{(\tau_{j-1} + t) \wedge \tau_{j}} $.
\end{lemma}

\begin{proof}
	Let us fix $x_1,x_2,\dotsc$ be a sequence of points in the grid.
	Let us define 
	\begin{equation*}
		Q(x;x_1,x_2,\dotsc)
		:=\Prob_x\Big(X_{\tau_1}=x_1,X_{\tau_2}=x_2,\dotsc\Big)
	\end{equation*}
	By the strong Markov property, 
	\begin{multline*}
		\Prob_x\Big(X_{\tau_{i+1}}=x_{i+1},X_{\tau_{i+2}}=x_{i+2},\dotsc
		\Big| \mathcal{F}_{\tau_i}\Big)
		=
		\Prob_{X_{\tau_i}}\Big(X_{\tau_{1}}=x_{i+1},,X_{\tau_{i+2}}=x_{i+2}\Big)
		\\
		=Q(X_{\tau_i};x_{i+1},x_{i+2},\dotsc).
	\end{multline*}
	Using the strong Markov property twice, first by 
	conditioning first with respect to $\mathcal{F}_{\tau_{j+1}}$
	and then with respect to $\mathcal{F}_{\tau_j}$, 
	\begin{multline*}
		\Esp_x\Big(
		\indic{X_{\tau_1}=x_1,
			\dotsc,
			X_{\tau_j}=x_j}
		F(\sandwitchprocessusage{X}{j})
		\indic{X_{\tau_{j}}=x_{j}, X_{\tau_{j+1}=x_{j+1},\dotsc}}
		\Big)
		\\
		=
		\Esp_x\Big(
		\indic{X_{\tau_1}=x_1,
			\dotsc,
			X_{\tau_j}=x_j}
		F(\sandwitchprocessusage{X}{j})
		\indicB{X_{\tau_{j}}=x_{j}}
		Q(X_{\tau_{j}})\Big) 
		\\
		=\Esp_x\Big(
		\indic{X_{\tau_1}=x_1,
			\dotsc,
			X_{\tau_{j-1}}=x_{j-1}}
		\sandwitchprocessusage{X}{j}
		\Big)Q(x_{j})
		=
		\Esp_x\Big(
		\indic{X_{\tau_1}=x_1,
			\dotsc,
			X_{\tau_j}=x_{j-1}}
		\Big)R(x_{j-1},x_{j})Q(x_{j})
	\end{multline*}
	with $R(x,y):=\Prob_x\Bigbraces{F((X^{0,\tau_1}_t)_{t\ge  0})\indic{X_{\tau_1}=y}}$.
	Therefore, using the definition of the conditional expectation, 
	\begin{multline*}
		\Esp_x\Big(
		F(\sandwitchprocessusage{X}{j})
		\Big|
		X_{\tau_1}=x_1,X_{\tau_2}=x_2,\dotsc\Big)
		\\
		=\frac{R(x_{j-1},x_{j})}{\Prob_{x_{j-1}}\Big(X_{\tau_1}=x_{j}\Big)}
		=
		\Esp_{x_{j-1}}\Big(F((X^{0,\tau_1}_t)_{t\geq 0}\Big|X_{\tau_1}=x_{j}\Big).
	\end{multline*}
	This is sufficient to prove the result.
\end{proof}

\paragraph{Acknowledgments}
The PhD thesis of  A.~Anagnostakis is supported by a scholarship from
the Grand-Est Region (France). 
